\documentclass{article}
\usepackage[utf8]{inputenc}
\usepackage{graphicx} 
\usepackage{amsthm}
\usepackage{amssymb}
\usepackage[UKenglish]{babel}
\usepackage{colonequals}
\usepackage{mathtools}

\usepackage{amsmath}	
\usepackage{amsfonts}   
\usepackage{enumitem}   

\usepackage{amsthm}
\usepackage{fancyhdr}   
\usepackage{MnSymbol}   
\usepackage[colorlinks]{hyperref}
\usepackage{tikz}
\usepackage{tikz-cd}    
\usetikzlibrary{calc}
\usepackage{comment}

\newtheorem{lemma}{Lemma}[section]
\newtheorem{thm}[lemma]{Theorem}
\newtheorem{cor}[lemma]{Corollary}
\newtheorem{question}[lemma]{Question}

\newtheorem*{lemma*}{Lemma}
\newtheorem*{lem*}{Lemma}
\newtheorem*{prop*}{Proposition}
\newtheorem*{thm*}{Theorem}

\theoremstyle{definition}

\newtheorem{defi}[lemma]{Definition}

\newtheorem{rem}[lemma]{Remark}
\newtheorem{example}[lemma]{Example}
\newtheorem*{rem*}{Remark}

\newenvironment{proofof}[1]{\begin{proof}[Proof of #1]}{\end{proof}}

\newcommand{\IR}{\mathbb{R}}

\newcommand{\IZ}{\mathbb{Z}}
\newcommand{\IN}{\mathbb{N}}
\newcommand{\IS}{\mathbb{S}}
\newcommand{\IF}{\mathbb{F}}
\newcommand{\IC}{\mathbb{C}}

\newcommand{\iso}{\cong}

\newcommand{\id}{\operatorname{id}}

\newcommand{\catC}{\mathcal{C}}
\newcommand{\Hom}{\operatorname{Hom}}

\newcommand{\Pol}{\operatorname{Pol}}
\newcommand{\proj}{\operatorname{proj}}
\newcommand{\conv}{\operatorname{conv}}

\newcommand{\CSP}{\operatorname{CSP}}

\renewcommand{\subset}{\subseteq}

\renewcommand{\epsilon}{\varepsilon}
\newcommand{\grp}{\underline{G}}

\newcommand{\de}[1]{\textbf{#1}}

\newcommand{\homProj}{\operatorname{HomProj}}
\newcommand{\Homsc}{\operatorname{Hom}^\text{SC}}
\newcommand{\NOT}{\operatorname{NOT}}
\newcommand{\dOR}{\operatorname{3-OR}}
\newcommand{\yvect}{\vec{y}}

\newcommand{\complexA}{A}
\newcommand{\complexB}{B}
\newcommand{\complexC}{C}
\newcommand{\subsetA}{\alpha}
\newcommand{\subsetAa}{\alpha'}
\newcommand{\partialmap}{\rho}

\newcommand{\facesA}{F_\complexA}
\newcommand{\facesB}{F_\complexB}
\newcommand{\face}{F}   
\newcommand{\dimension}{d}
\newcommand{\morph}{f}
\newcommand{\pathcomplex}{\mathcal{P}}    
\newcommand{\cyclecomplex}{\mathcal{C}}    
\newcommand{\georeal}[1]{\operatorname{geom}(#1)}   
\newcommand{\topspaceA}{X}
\newcommand{\topspaceB}{Y}
\newcommand{\universal}{Hom-universal}  
\newcommand{\contractible}{Hom-contractible}  
\newcommand{\siggers}{S}    
\newcommand{\cyclic}{c}    

\newcommand{\dSAT}{\operatorname{3-SAT}} 

\newcommand{\struct}[1]{\underline{#1}}
\newcommand{\structA}{{\struct{A}}}
\newcommand{\signature}{\sigma}
\newcommand{\setstructA}{A}
\newcommand{\structB}{{\struct{B}}}
\newcommand{\subsetB}{\beta}
\newcommand{\setstructB}{B}
\newcommand{\structC}{\struct{C}}
\newcommand{\setstructC}{C}
\newcommand{\subsetC}{\gamma}
\newcommand{\structD}{\struct{D}}
\newcommand{\Relationsymbol}{R}
\newcommand{\clone}{C}
\newcommand{\clonefunctiona}{c}
\newcommand{\clonefunctionb}{d}
\newcommand{\clonefunctionc}{e}
\newcommand{\clonefunctions}{s}
\newcommand{\Polidem}{\operatorname{Pol}_{\text{idem}}} 
\newcommand{\Proj}{\operatorname{proj}} 

\newcommand{\facerel}{F} 
\newcommand{\fix}{_{\textup{fix}}}
\newcommand{\core}{\operatorname{core}}
\newcommand{\cloneD}{D} 
\newcommand{\subdivision}[1]{\operatorname{GS}(#1)} 
\newcommand{\subdivisionn}[1]{\operatorname{GS}^n(#1)} 
\newcommand{\auto}{\gamma} 

\title{A Dichotomy for \\ Finite Abstract Simplicial Complexes}
\author{Sebastian Meyer\footnote{Funding statement: Funded by the European Union (ERC, POCOCOP, 101071674). Views and opinions expressed are however those of the author(s) only and do not necessarily reflect those of the European Union or the European Research Council Executive Agency. Neither the European Union nor the granting authority can be held responsible for them. }}
\date{\today}

\begin{document}

\maketitle
\begin{abstract}
    Given two finite abstract simplicial complexes $\complexA$ and $\complexB$, one can define a new simplicial complex on the set of simplicial maps from $\complexA$ to $\complexB$. After adding two technicalities, we call this complex $\Homsc_{\subsetA,\partialmap}(\complexA, \complexB)$.


    We prove the following dichotomy: For a fixed finite abstract simplicial complex $\complexB$, either $\Homsc_{\subsetA,\partialmap}(\complexA, \complexB)$ is always a disjoint union of contractible spaces or every finite CW-complex can be obtained up to a homotopy equivalence as $\Homsc_{\subsetA,\partialmap}(\complexA, \complexB)$ by choosing $\complexA$ in a right way. 
    
    We furthermore show that the first case is equivalent to the existence of a nontrivial social choice function and that in this case, the space itself is homotopy equivalent to a discrete set.
    
    Secondly, we give a generalization to finite relational structures and show that this dichotomy coincides with a complexity theoretic dichotomy for constraint satisfaction problems, namely in the first case, the problem is in P and in the second case NP-complete. 
    This generalizes a result from \cite{BooleanCaseLIPIcs}.
    %
\end{abstract}
\tableofcontents
\section{Introduction}

For simplicial complexes and other topological spaces, one can notice that more symmetries often imply a simpler homotopy type. The most common example of this fact is the theorem that every connected topological group has a commutative fundamental group. If we add the assumption that the addition of the group should also be a simplicial map, we even get that the underlying space is contractible.

This fact extends to a dichotomy for abstract simplicial complexes. If a complex $\complexB$ has a simplicial map from a power of $\complexB$ to $\complexB$ which satisfies some weak symmetry conditions, then $\complexB$ is a disjoint union of contractible spaces. Moreover all simplicial complexes that can be constructed from this one as homomorphisms into $\complexB$ will also be disjoint unions of contractible spaces. This type of symmetry maps are called Taylor polymorphisms and they always exist for groups. On the other side of the dichotomy, we get that a complex without this type of symmetry can construct this way actually all finite simplicial complexes up to a homotopy equivalence. We give the precise statement in Section~\ref{SectionMainTheorem}.

Besides this purely topological statement, this dichotomy is also closely connected with two other areas of mathematics: social choice functions and constraint satisfaction problems.



\subsection{Social Choice Functions}
\label{SectionSocialChoiceIntro}
Social choice functions are essentially the same as polymorphisms. Thus the dichotomy is also related to the existence of social choice functions. We give a more concrete example:

Consider a big conference with a set of participants and a set $A$ of possible schedules. Every two people with the same schedule will have the chance to be in contact at this conference, but there are also some shared events, so also people with different schedules will have the chance to get in contact. We define a subset of $F\subset A$ to be a face if there is an event on which all people from schedules on $F$ meet at the same time. It is easy to see that this defines a simplicial complex on $A$.

\begin{example} \label{ExampleMotivationConference}
    Consider a conference with a big and a small seminar room. The conference has three days and there are three possible schedules. In each day, two schedules visit the big seminar room and one is in the small one such that every schedule contains two days in the big room and one day in the small one.
    The simplicial complex drawn from this conference is an empty triangle as every two schedules have a common day but there is no day on which all three schedules are together. This complex is not contractible.
\end{example}

Now consider any number $k$ of organizers to this conference. Each of them proposes an assignment of participants to schedules. We assume furthermore that the conference is big enough so every possible combination of $k$ assignments on $A$ actually occurs. Now, it would be nice to have a final assignment, satisfying all of the conditions
\begin{enumerate}
    \item Unanimity: If all $k$ proposals agree for a person, then this person should get this schedule.
    \item Mapping property: If two persons get the same proposals, they should get the same schedule. In other words, the final schedule should only depend on the proposals, not on the participant.
    \item Homomorphism property: If a set of participants meet each other at the same time in each proposal, then they should do so in the final assignment. That represents the fact that they likely want to work together in a group and this is known and respected by the organizers.
\end{enumerate}
and some conditions like the following:
\begin{enumerate}[resume]
    \item Anonymity: The schedules do not change when the order of the proposals is changed.
    \item Majority: If more than halve of the proposals agree on the schedule for a person, this person should get that schedule.
    \item Near unanimity: If $k\ge 3$ and all but one of the proposals agree on the schedule for a person, this person should get the schedule of the majority. (This is implied by majority.)
\end{enumerate}
The first three condition can be easily satisfied by declaring one of the organizers to be the dictator and just take her or his proposal.
Adding one of the next three conditions makes it much harder to find a solution.

We show that the first three conditions together with any Taylor condition, that is essentially any condition which is not satisfied by the dictator algorithm, cannot be satisfied if the simplicial complex of the conference is not contractible. Therefore, there is no choice algorithm in Example~\ref{ExampleMotivationConference} to choose a final assignment which satisfies the first three and any of the last three conditions.

For a majority on three organizers, this was already shown in \cite[Theorem 7.7]{Taylor77}. For anonymity, it was shown in \cite{GeneralizedMeans} and \cite{SocialChoice} that the space has to be contractible under the additional assumption that there is such a choice function for every $k$.

\subsection{Constraint Satisfaction Problems}
A constraint satisfaction problem (CSP) is a classical computer problem. It is described by a template that is a relational structure. As input, it get a primitive positive formula which can include some free variables and some existentially quantified variables. As output, it should return a possible assignment for the free variables in the template that satisfies the formula respectively return that there is no such assignment.


CSPs include solving a HORNSAT-formula or a SAT-formula. Also determining wether a graph is three colorable is a CSP where the variable set is the set of vertices and the template the set of colors. Furthermore, it is also a CSP to determine if any polynomial can be solved over a fixed ring or for a fixed structure to check, if any finite second structure has a homomorphism into the first one. 
In fact every CSP can be formalized in the following form: 
Given a signature $\signature$ and a template $\signature$-structure $\structB$, determine for any given finite input $\signature$-structure $\structA$, if there is a $\sigma$-homomorphism $\structA\to \structB$.

It is easy to see that if the template is finite, then the CSP can be solved in exponential time and furthermore in NP because if there is a correct assignment, one can just guess the values for all variables and verify it easily. Moreover, some CSPs such as solving HORNSAT formulas, are known to be in P. It was shown in \cite{Bulatov} and \cite{Zhuk}, that there is a complexity theoretic dichotomy for finite CSPs:

\begin{enumerate}
    \item Either the structure $\structB$ has a Taylor polymorphism. Then, the CSP of $\structB$ can be solved in polynomial time.
    \item Or the structure $\structB$ has no Taylor polymorphism. Then, the CSP of $\structB$ is NP-complete.
\end{enumerate}

We can extend the dichotomy for simplicial complexes given in the beginning to arbitrary structures. This gives an equivalent topological characterization of the complexity theoretic dichotomy. To define this generalization, fix any structure $\structB$. Consider an instance of its CSP in the version described above. 
Look at the set of all solutions of the free variables. This set carries naturally the structure of a finite simplicial complex. Now, we get the dichotomy for finite CSPs:
\begin{enumerate}
    \item
    Either, for every instance, the resulting space is a contractible. In this case, $\structB$ has a Taylor polymorphism and its CSP is in P-time.
    \item 
    Or for every finite simplicial complex, there is an instance such that the solution space is homotopy equivalent to the chosen complex. In this case, $\structB$ has no Taylor polymorphism and its CSP is NP-complete.
\end{enumerate}
We make this more formal in Theorem~\ref{TheoremDualityForStructures}.
The idea of this dichotomy is based on \cite{BooleanCaseLIPIcs}, in which a slightly different topological dichotomy was shown in the special case of boolean CSPs.

\section{General Framework}
This section presents the general topological framework which is needed for a good understanding of the main theorem. 

\subsection{Abstract Simplicial Complexes}

\begin{defi}
    We define an \de{(abstract) simplicial complex} to be a set $\complexA$ together with a subset $\facesA$ of the power set $\powerset(\complexA)$ of $\complexA$ called \de{faces} which contains only finite sets, all subsets of size $1$ and the empty set and is closed under taking subsets.

    An abstract simplicial complex $(\complexA,\facesA)$ is called \de{finite} if $\complexA$ is finite. It is called \de{finite dimensional} if there is a natural number $\dimension$ such that $\facesA$ contains only sets with size at most $\dimension+1$. The smallest such $\dimension$ is called \de{dimension} of $(\complexA,\facesA)$ and is written $\dim(\complexA)$.
\end{defi}

We will restrict our attention later to finite simplicial complexes. Every finite simplicial complex $(\complexA,\facesA)$ is finite dimensional with dimension at most $|\complexA|-1$.
The following two examples will play a role later.
\begin{example} \label{ExampleCyclepath}
    The \de{path of length $n$} is the simplicial complex $\pathcomplex_n$ on the set $\{0,1,...,n\}$ with the faces $\emptyset$, $\{\{i\}\mid 0\le i\le n\}$ and $\{\{i,i+1\}\mid 0\le i\le n-1\}$.

    The \de{cycle of length $n$} is the simplicial complex $\cyclecomplex_n$ on the set $\{1,...,n\}$ with the faces $\emptyset$, $\{\{i\}\mid 1\le i\le n\}$ and $\{\{i,i+1\}\mid 1\le i\le n-1\}\cup\{\{n,1\}\}$.
\end{example}

\begin{defi}[Products, Homomorphisms]
    Let $(\complexA,\facesA)$ and $(\complexB,\facesB)$ be finite simplicial complexes. 
    
    Their \de{product} is defined as the set $\complexA \times \complexB$ where a set $\{(x_i,y_i)\mid i\in I\}\subset \complexA \times \complexB$ is a face if and only if both $\{x_i\mid i\in I\}\subset \complexA$ and $\{y_i\mid i\in I\}\subset \complexB$ are faces. The product is written $\complexA \times \complexB$. 


    A \de{simplicial homomorphism} from $(\complexA,\facesA)$ to $(\complexB,\facesB)$ is a map $\morph\colon  \complexA \to \complexB$ such that for all faces in $\facesA$, their image is in $\facesB$. We call the set of all these maps $\Hom(\complexA,\complexB)$. 
\end{defi}
The product satisfies the universal property of the product. If $\complexA$ and $\complexB$ are both finite respectively finite dimensional then so is their product. Note that the dimension of the product is only bounded by $(\dim(\complexA)+1)(\dim (\complexB)+1)-1$ which is sharp for example for $\facesA=\powerset(\complexA)$ and $\facesB=\powerset(\complexB)$. Here and in the rest of this paper, we denote the power set of a set $X$ by $\powerset(X)$ and the size of the set by $|X|$. 
The geometric realization, which we introduce later, will be denoted as $\georeal{X}$ and not by $|X|$ to not raise any confusions.

We can make the space of simplicial homomorphisms again a simplicial complex. We moreover have even two meaningful definitions: 
\begin{defi}
    Let $(\complexA,\facesA)$ and $(\complexB,\facesB)$ be finite simplicial complexes. We identify the set of all maps from $\complexA$ to $\complexB$ with the product $\complexB^{|\complexA|}$. 

    We define a simplicial set on the set $\Hom(\complexA,\complexB)$ by declaring a subset $\subsetA \subset \Hom(\complexA,\complexB)$ to be a face, if the corresponding subset of $\complexB^{|\complexA|}$ is a face. We name this complex again $\Hom(\complexA,\complexB)$. This construction is called \de{internal hom-functor} (on objects). 

    To define the \de{simplicial-hom-complex} $\Homsc(\complexA,\complexB)$, we define a simplicial structure on $\Hom(\complexA,\complexB)$ as follows:
    We identify again $\Hom(\complexA,\complexB)$ with tuples in $\complexB^{|\complexA|}$. We say a set
    $\{ (x_{i,1},...,x_{i,|\complexA|}) \mid i \in \{1,...,k\}\}\subset B^{|A|}$ of size $k$ is a face if for all maps $j\colon \{1,...,|\complexA|\}\to \{1,...,k\}$, the tuple $(x_{j(1),1},...,x_{j(|\complexA|),|\complexA|})$ corresponds to a homomorphism.
\end{defi}

\begin{example}\label{Example1}
    Let 
    \begin{align*}
        (\complexA,\facesA)&=(\{1,2\}, \{\emptyset, \{1\},\{2\},\{1,2\}\})\iso \pathcomplex_1 \text{ and} \\ 
        (\complexB,\facesB)&=(\{x,y,z\}, \{\emptyset, \{x\},\{y\},\{z\},\{x,y\},\{y,z\}\})\iso \pathcomplex_2,
    \end{align*}
    see Figure~\ref{FigureExample1} for the drawing. For the map $1\mapsto a$ and $2\mapsto b$, we simply write $ab$. With this notation, we get
    \begin{align*}
        \Hom( \complexA, \complexB) &=(\{xx,xy,yx,yy,yz,zy,zz\}, 
        \\&\phantom{{}={(}}\phantom{{}\cup{}} \powerset(\{xx,yy,xy,yx\}) \cup \powerset(\{yy,zz,yz,zy\})
        \\& \phantom{{}={(}}\cup \powerset(\{xy,yy,yz\})\cup \powerset(\{yx,yy,zy\})
        )\text{ and} \\ 
        \Homsc( \complexA, \complexB) &=(\{xx,xy,yx,yy,yz,zy,zz\}, 
        \\&\phantom{{}={(}}\phantom{{}\cup{}} \powerset(\{xx,yy,xy,yx\}) \cup \powerset(\{yy,zz,yz,zy\})
        \\& \phantom{{}={(}} \cup \powerset(\{xy,yy,zy\})\cup \powerset(\{yx,yy,yz\})
        ).
    \end{align*}
    They clearly have the same base set but differ on the faces.
    \begin{figure}
        \centering
        \begin{tikzpicture}
            \draw (0,0) node[anchor=north] {$1$}-- (2,0) node[anchor=north] {$2$};
            \node at (1,-1) {$\complexA$};
            \draw (3,0) node[anchor=north] {$x$} -- (4,2) node[anchor=south] {$y$} -- (5,0) node[anchor=north] {$z$};
            \node at (4,-1) {$\complexB$};
            \fill [fill=gray!20!white] (6,2) -- (7,2) -- (8,1) -- (8,0) -- (7,0) -- (6,1) -- cycle;
            \draw (6,2) -- (7,2) -- (7,0) -- (8,0) -- (8,1) -- (6,1) -- cycle;
            \draw (7,2)--(6,1)--(7,0)--(8,1)--cycle
                (6,2)--(8,0);    
            \draw (6,2) node [fill=white]{$xx$} 
                (6,1) node[fill=white] {$yx$} 
                (7,2) node[fill=white] {$xy$} 
                (7,1) node[fill=white] {$yy$} 
                (7,0) node[fill=white] {$zy$} 
                (8,1) node[fill=white] {$yz$} 
                (8,0) node[fill=white] {$zz$};
            \node at (7,-1) {$\Hom(\complexA,\complexB)$};
            \fill[fill=gray!20!white] (9,1) to [bend right] (11,1) -- cycle
                 (10,2) to [bend left=40] (10,0) -- cycle;   
            \draw[fill=gray!20!white] (9,2) -- (10,2) -- (10,0) -- (11,0) -- (11,1) -- (9,1) -- cycle;
            \draw (10,2)--(9,1) 
                (10,0)--(11,1) 
                (9,2)--(11,0)
                (9,1) to [bend right] (11,1)
                (10,2) to [bend left=40] (10,0);    
            \draw (9,2) node [fill=white]{$xx$} 
                (9,1) node[fill=white] {$yx$} 
                (10,2) node[fill=white] {$xy$} 
                (10,1) node[fill=white] {$yy$} 
                (10,0) node[fill=white] {$zy$} 
                (11,1) node[fill=white] {$yz$} 
                (11,0) node[fill=white] {$zz$};
            \node at (10,-1) {$\Homsc(\complexA,\complexB)$};
        \end{tikzpicture}
        \caption{A graphical visualization of the complexes from Example~\ref{Example1}. Non-degenerate faces with exactly two elements are drawn as lines. Higher dimensional faces are highlighted in gray. Note that some of them overlap which is not visible in the picture.}
        \label{FigureExample1}
    \end{figure}
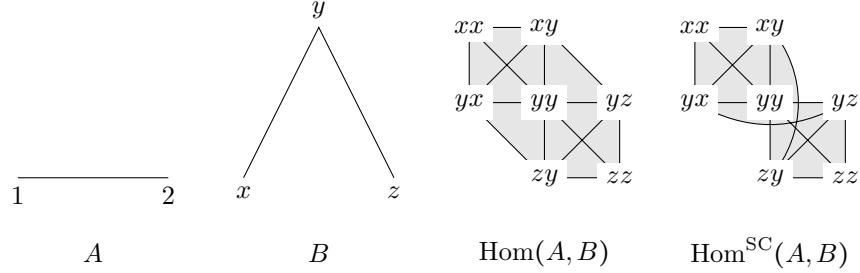
\end{example}

Both constructions are functorial in the sense that if there are simplicial maps $\complexA' \to \complexA$ and $\complexB \to \complexB'$, then thy induce simplicial maps $\Homsc(\complexA,\complexB) \to \Homsc(\complexA',\complexB')$ and $\Hom(\complexA,\complexB) \to \Hom(\complexA',\complexB')$ by composition.

On the side of differences, we have that only
the internal Hom functor has a left adjoint tensor product, while 
the composition of maps induces a simplicial morphism from $\Homsc(\complexB, \complexC)\times \Homsc(\complexA,\complexB)$ to $\Homsc(\complexA,\complexC)$. It thus behaves more like a construction method and not like the internal Hom functor. 
When generalizing the simplicial-Hom-complex to other categories $\catC$, one should consider it as a (bi-)functor from any $\catC$ into simplicial complexes while the internal Hom functor would again map into $\catC$.

We want to enlarge the possible constructions in two ways by introducing fixed points and some kind of projection: 
\begin{defi}
    Let $(\complexA,\facesA)$ and $(\complexB,\facesB)$ be two simplicial complexes. Let $\subsetA$ and $\subsetAa$ be two subsets of $\complexA$. Let $\partialmap\colon \subsetAa \to \complexB$ any map. 
    
    Then, define $\Hom_{\subsetA,\partialmap}$ as the full subcomplex of $\complexB^{|\subsetA|}$ consisting of all points which correspond to a function $\subsetA \to \complexB$ which has an extension to $\complexA$ that also extends extends $\partialmap$ and is a homomorphism.

    Define $\Homsc_{\subsetA,\partialmap}$ as a simplicial set: The set is the subset of $\complexB^{|\subsetA|}$ consisting of all points which correspond to a function $\subsetA \to \complexB$ which can be extended to a simplicial homomorphism $\complexA \to \complexB$ that also extends $\partialmap$. A finite subset $\{ (x_{i,1},...,x_{i,|\complexA|}) \mid i \in \{1,...,k\}\}$  of $\Hom_{\subsetA,\partialmap}$ is a face if for all maps $j\colon \{1,...,|\complexA|\}\to \{1,...,k\}$, the tuple $(x_{j(1),1},...,x_{j(|\complexA|),|\complexA|})$ is again in $\Homsc_{\subsetA,\partialmap}$.

    If $\subsetAa=\emptyset$, we write $\partialmap=\emptyset$.
\end{defi}

This definition turns out to be less intuitive as one might think:
\begin{example} \label{Example2}
    Consider
    \begin{align*}
        (\complexA,\facesA) &=\pathcomplex_2= (\{0,1,2\}, \{\{\},\{0\},\{1\},\{2\},\{0,1\},\{1,2\}\}) 
        \\
        (\complexB,\facesB) &=\cyclecomplex_5= (\{1,2,3,4,5\}, ...) 
        \\
        \subsetA &= \{2\} \text{ and }
        \partialmap \colon \{0\} \to \complexB , 0\mapsto 1.
    \end{align*}
    Then, we get that the equivalence classes of the homomorphisms $(0,1,2)\mapsto(1,2,3)$ and $(0,1,2)\mapsto(1,5,4)$ are connected in both $\Hom_{\subsetA,\partialmap}(\complexA,\complexB)$ and $\Homsc_{\subsetA,\partialmap}(\complexA,\complexB)$, but they are not connected in $\Hom_{\complexA,\partialmap}(\complexA,\complexB)$ nor in $\Homsc_{\complexA,\partialmap}(\complexA,\complexB)$, see Figure~\ref{FigureExample2}. It turns out that the surjective maps are no quotients.
    
    \begin{figure}
        \centering
        \begin{tikzpicture}
            \draw [fill=gray!20!white] (-1,0) -- (1,0) -- (3,2) -- (1,2) -- cycle;
            \draw (0,1) -- (2,1)
                (0,0) -- (2,2)
                (0,0) -- (0,1)
                (1,0) -- (1,2)
                (2,1) -- (2,2)
                (1,2) -- (2,1)
                (0,1) -- (1,0);
            \draw (-1,-2) to [bend right] (3,-2) -- cycle;
            \draw (1,1) node [fill=white]{$111$} 
                (1,2) node[fill=white] {$121$} 
                (1,0) node[fill=white] {$151$} 
                (2,1) node[fill=white] {$112$} 
                (2,2) node[fill=white] {$122$} 
                (3,2) node[fill=white] {$123$} 
                (0,1) node[fill=white] {$115$} 
                (0,0) node[fill=white] {$155$} 
                (-1,0) node[fill=white] {$154$};
            \node at (1,-1) {$\Hom_{\complexA,\partialmap}(\complexA,\complexB)$};
            \draw (1,-2) node [fill=white]{$1$} 
                (2,-2) node[fill=white] {$2$} 
                (3,-2) node[fill=white] {$3$} 
                (0,-2) node[fill=white] {$5$} 
                (-1,-2) node[fill=white] {$4$};
            \node at (1,-3) {$\Hom_{\subsetA,\partialmap}(\complexA,\complexB)$};
            \fill [fill=gray!20!white] (5,0) to [bend right] (7,0) -- cycle;
            \fill [fill=gray!20!white] (6,1) to [bend right] (8,1) -- cycle;
            \fill [fill=gray!20!white] (7,2) to [bend right] (9,2) -- cycle;
            \fill [fill=gray!20!white] (7,0) to [bend right=40] (7,2) -- cycle;  
            \fill [fill=gray!20!white] (6,1) --(8,1) -- (8,2)--(7,2)--(7,0)--(6,0)--cycle;
            \draw (5,0) to [bend right] (7,0) -- cycle;
            \draw (6,1) to [bend right] (8,1) -- cycle;
            \draw (7,2) to [bend right] (9,2) -- cycle;
            \draw (7,0) to [bend right=40] (7,2) -- cycle;  
            \draw (6,1) -- (6,0)
                (8,2)--(8,1)
                (8,2)--(6,0)
                (7,2)--(8,1)
                (6,1)--(7,0);
            \draw (7,1) node [fill=white]{$111$} 
                (7,2) node[fill=white] {$121$} 
                (7,0) node[fill=white] {$151$} 
                (8,1) node[fill=white] {$112$} 
                (8,2) node[fill=white] {$122$} 
                (9,2) node[fill=white] {$123$} 
                (6,1) node[fill=white] {$115$} 
                (6,0) node[fill=white] {$155$} 
                (5,0) node[fill=white] {$154$};
            \node at (6,-1) {$\Homsc_{\complexA,\partialmap}(\complexA,\complexB)$};
            \draw [fill=gray!20!white] (5,-2) to [bend right] (9,-2) -- cycle;
            \draw (5,-2) to [bend right] (7,-2)
                (5,-2) to [bend right] (8,-2)
                (6,-2) to [bend right] (8,-2)
                (6,-2) to [bend right] (9,-2)
                (7,-2) to [bend right] (9,-2);
            \draw (7,-2) node [fill=white]{$1$} 
                (8,-2) node[fill=white] {$2$} 
                (9,-2) node[fill=white] {$3$} 
                (6,-2) node[fill=white] {$5$} 
                (5,-2) node[fill=white] {$4$};
            \node at (6,-3) {$\Homsc_{\subsetA,\partialmap}(\complexA,\complexB)$};
        \end{tikzpicture}
        \caption{The simplicial complexes from Example~\ref{Example2}.}
        \label{FigureExample2}
    \end{figure}
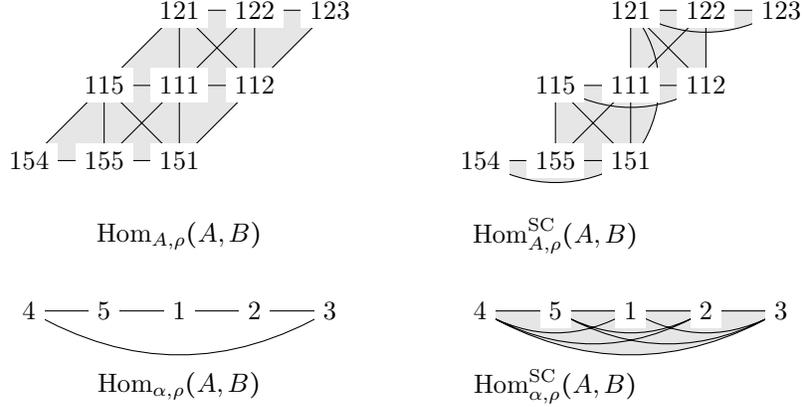
\end{example}
\subsection{Geometric Definitions}
There is a forgetful functor from abstract simplicial complexes to topological spaces called geometrical realization. We will use this connection to define when we call an abstract simplicial complex contractible and to define when complexes are homotopy equivalent.

We define this construction only for a finite complex.
\begin{defi} \label{DefinitionGeometricRealization}
    The geometric realization of a finite abstract simplicial complex $(\complexA,\facesA)$ is the subset of $\IR^{|\complexA|}$ given by
    $$
        \bigcup_{\face \in \facesA} \conv \{e_i \mid i \in \face\} 
    $$
    where $e_i\in \IR^{|\complexA|}$ is the unit basis vector corresponding to $i\in \complexA$ and $\conv$ denotes the convex hull. The $n$-skeleton is the subset
    $$
        \bigcup_{\face \in \facesA, |\face|\le n+1} \conv \{e_i \mid i \in \face\} 
    $$
    which consists only of the at most $n$-dimensional faces. The geometric realization of $\complexA$ will be denoted $\georeal{\complexA}$.
\end{defi}
It can be easily seen that for a finite complex $\complexA$, the realization $\georeal{\complexA}$ will always be a compact Hausdorff space and also a CW-complex.

Moreover, the geometric realization is a functor by extending maps from the vertices linearly. It does not preserves products as the path $\pathcomplex_1$ is realized by a compact interval and the product $\pathcomplex_1 \times \pathcomplex_1$ is realized by a tetrahedron. However, the product of the topological spaces is a subset of the simplicial product so a simplicial map $\complexA^n \to \complexB$ induces a continuous map $\georeal{\complexA}^n \to \georeal{\complexB}$.

Moreover, $\georeal{\complexA}^n$ and $\georeal{\complexA^n}$ are homotopy equivalent, which is the next definition:
\begin{defi}
    Let $\topspaceA$ and $\topspaceB$ be two topological spaces. We call two continuous maps  $f,g\colon \topspaceA \to \topspaceB$ \de{homotopic}, if there is a third continuous map $F\colon \topspaceA \times [0,1] \to \topspaceB$ such that $f(x)=F(x,0)$ and $g(x)=F(x,1)$ for all $x\in \topspaceA$.

    The spaces $\topspaceA$ and $\topspaceB$ are \de{homotopy equivalent}, if there are two continuous maps $f\colon \topspaceA \to \topspaceB$ and $g\colon \topspaceB \to \topspaceA$ such that $g\circ f$ is homotopic to the identity on $\topspaceA$ and $f\circ g$ is homotopic to the identity on $\topspaceB$.

    The space $\topspaceA$ is called \de{contractible}, if it is homotopy equivalent to a single point.
\end{defi}
We call abstract simplicial complexes \de{homotopy equivalent} or \de{contractible} if their geometric realizations have the respective property.

There are two important cases how the spaces $\Homsc(\complexA,\complexB)$ can behave up to homotopy:
\begin{defi} \label{DefinitionContractible}
    Let $\complexB$ be a finite abstract simplicial complex. 
    
    We call $\complexB$ \de{\contractible{}}, if for every finite abstract simplicial complex $\complexA$, 
    every subsets $\subsetA, \subsetAa\subset \complexA$ and every partially defined map $\partialmap\colon \subsetAa \to \complexB$, 
    the complex $\Homsc_{\subsetA,\partialmap}(\complexA,\complexB)$ has only contractible components.

    We call $\complexB$ \de{\universal}, if for every finite abstract simplicial complex $\complexC$, 
        there is a finite abstract simplicial complex $\complexA$, subsets $\subsetA, \subsetAa\subset \complexA$ and a map $\partialmap\colon \subsetAa \to \complexB$, such that the complex $\Homsc_{\subsetA,\partialmap}(\complexA,\complexB)$ is homotopy equivalent to $\complexC$.
\end{defi}

\section{Main Theorem}
\label{SectionMainTheorem}

We are now able to state our main theorem:

\begin{thm} \label{TheoremMain}
    Every finite abstract simplicial complex is either \contractible{} or \universal{}.

\end{thm}

There are many consequences of being \contractible{}. By a contraposition, we get some criteria for being \universal, which include the following.
\begin{thm} \label{TheoremExtra}
    Let $\complexB$ be a finite abstract simplicial complex. If
    \begin{enumerate}
        \item $\complexB$ has a non-contractible connected component or
        \item there are any finite $(\complexA,\subsetA,\partialmap)$ such that $\Hom_{\subsetA,\partialmap}(\complexA,\complexB)$ has a non-contractible connected component,
    \end{enumerate}
	then $\complexB$ is \universal{}.
\end{thm}
This property can also be described in the existence of some simplicial homomorphisms. We want to give them already now, even though some of the definitions only appear later.
\begin{thm} \label{TheoremExtra2}
	The following are equivalent for a simplicial complex $\complexB$.
    \begin{enumerate}
    	\item The complex $\complexB$ is \universal{}.
        \item The complex $\complexB$ is not \contractible{}.
        \item The complex $\complexB$ has an idempotent Taylor polymorphism.
        \item \label{NrTheoremExtraSiggers}
        There is an idempotent Siggers polymorphism \cite{Siggers}, that is a simplicial homomorphism $\siggers\colon \complexB^6\to \complexB$ such that 
        \begin{align*}
            \forall x \in \complexB: &&\siggers(x,x,x,x,x,x)&=x \text{ and}\\
            \forall x,y,z \in \complexB: && \siggers(x,x,y,y,z,z)&=\siggers(z,y,x,z,y,x)
        \end{align*}
        hold.
        \item There are idempotent cyclic Taylor polymorphism, that are simplicial homomorphisms $\cyclic\colon \complexB^n\to \complexB$  for all prime numbers $n>|\complexB|$ such that 
        \begin{align*}
            \forall x \in \complexB: &&\cyclic(x,\dots ,x)&=x \text{ and}\\
            \forall x_0,x_1,...,x_{n-1} \in \complexB: && \cyclic(x_0,x_1,\dots,x_{n-1})&=\cyclic(x_1,\dots,x_{n-1},x_0)
        \end{align*}
        hold.
    \end{enumerate}
\end{thm}
There is even a complexity theoretic characterization:
\begin{thm} \label{TheoremExtra3}
        If $\complexB$ is \universal, then the decision problem given a finite $(\complexA,\subsetAa,\partialmap)$, determine whether $\Hom_{\complexA,\partialmap}(\complexA,\complexB)$ is nonempty, is NP-complete. Otherwise it is in P.
\end{thm}

Note that it is in general easy to prove that a space is not \contractible{} and thus \universal{}, because one needs to give only one example. To prove that 
a space is \contractible{}, one can use characterization \ref{NrTheoremExtraSiggers} of Theorem \ref{TheoremExtra2}.

The idea to this theorem is based on \cite[Theorem 10]{BooleanCase} in which a similar theorem for boolean constraint satisfaction problems was shown, which we generalize in Theorem~\ref{TheoremDualityForStructures} and on \cite{TopologyAndAdjunction} in which a similar theory was developed independently.

The proof of Theorem~\ref{TheoremMain} 
has the following key ideas:
We look at the idempotent polymorphisms of the abstract simplicial complex. The polymorphisms are a generalization to the endomorphism monoid and we introduce them more precisely in Section~\ref{sectionPolymorphisms}. Moreover, we can view the finite abstract simplicial complex $\complexB$ as a relational structure. In that language, there are many known Theorems about polymorphisms which we will be able to use, including a dichotomy linked with Taylor polymorphisms.


After that, there are three more steps to take:
\begin{itemize}
    \item Show that if $\complexB$ has a Taylor polymorphism, then also $\Hom_{\subsetA,\partialmap}(\complexA,\complexB)$ and $\Homsc_{\subsetA,\partialmap}(\complexA,\complexB)$ have Taylor polymorphisms. This is a simple fact shown in Corollary~\ref{KorollarClonhomozuHom}.
    \item Show that if $\complexB$ has a Taylor polymorphism (or anything equivalent) then every component of $\complexB$ is contractible. We show this in is Theorem~\ref{TheoremTaylorContractible}.
    \item Show that if $\complexB$ has no Taylor polymorphism, then it is \universal. This is done in Lemma~\ref{LemmaConstuctToUniversal} respectively Corollary~\ref{CorollaryNoTaylorUniversal}.
\end{itemize}

\begin{figure}
    \centering
    \begin{tikzpicture}
        \node [text width=5cm, align=center, draw] (left 2)
            at (0,-2) 
            {$\complexB$ has a Taylor \\ polymorphism} ;
        \node (test) [text width=5cm, align=center, draw] (left 3)
            at (0,-6) 
            {$\complexB$ has idempotent cyclic Taylor polymorphisms};
        \node [text width=3cm, align=center, draw, dashed] (left A)
            at (-1,-8)
            {The problem corresponding to $\complexB$ is in P-time};
        \node [text width=3cm, align=center, draw] (left B)
            at (-1,-4)
            {$\complexB$ is has an idempotent Siggers polymorphism};
        \node [text width=5cm, align=center, draw] (left 4)
            at (0,-10)
            {$B$, $\Hom_{\subsetA,\partialmap}(\complexA,\complexB)$ and $\Homsc_{\subsetA,\partialmap}(\complexA,\complexB)$ have idempotent cyclic Taylor polymorphisms};
        \node [text width=5cm, align=center, draw] (left 5)
            at (0,-12)
            {all components of $B$, $\Hom_{\subsetA,\partialmap}(\complexA,\complexB)$ and $\Homsc_{\subsetA,\partialmap}(\complexA,\complexB)$ are contractible};
        \node [text width=5cm, align=center, draw] (left 6)
            at (0,-14)
            {$B$ is \contractible{}};
        \node [text width=5cm, align=center, draw] (left 7)
            at (0,-16)
            {$B$ is not \universal};
    
        \node [text width=5cm, align=center, draw] (right 1) at (6,-2) {$\complexB$ has no Taylor \\ polymorphism} ;
        \node [text width=5cm, align=center, draw] (right 2) at (6,-6) {$\complexB$ pp-interprets $\dSAT$} ;
        \node [text width=5cm, align=center, draw] (right 3) at (6,-10) {$\complexB$ pp-interprets $\dSAT$ in one dimension} ;
        \node [text width=3cm, align=center, draw, dashed] (right B) at (7,-12) {$\dSAT$ is \universal{} \cite{BooleanCase}};
        \node [text width=3cm, align=center, draw, dashed] (right A) at (7,-8) {The problem corresponding to $\complexB$ is NP-complete};
        \node [text width=5cm, align=center, draw] (right 6) at (6,-14) {$B$ is \universal};

        \draw [<->, dashed] (left 2) -- (right 1);
        \draw [<->, dashed] (left 7) -| ($(right 6.south)+(-1,0)$);
        \draw [-Latex] ($(left 2.south)+(1,0)$) -- node [right]{\cite{CyclicTaylor,CyclicTaylorReprint}} ($(left 3.north)+(1,0)$);
        \draw [-Latex] ($(left 3.south)+(1,0)$) -- node [right]{Cor.~\ref{KorollarClonhomozuHom}}($(left 4.north)+(1,0)$);
        \draw [-Latex] ($(left 4.south)+(1,0)$) -- node [right]{Thm.~\ref{TheoremSimplicialCyclic}} ($(left 5.north)+(1,0)$);
        \draw [-Latex] ($(left 5.south)+(1,0)$) -- ($(left 6.north)+(1,0)$);
        \draw [-Latex] ($(left 6.south)+(1,0)$) -- ($(left 7.north)+(1,0)$);
        \draw [-Latex] ($(right 1.south)+(-1,0)$) -- ($(right 2.north)+(-1,0)$);
        \draw [-Latex] ($(right 2.south)+(-1,0)$) -- node [left, align=right]{Lem.~\ref{LemmaConsturctInterpret} \\ {\tiny (previously \tiny known)}}($(right 3.north)+(-1,0)$);
        \draw [-Latex] ($(right 3.south)+(-1,0)$) -- node [left, near end]{Lem.~\ref{LemmaConstuctToUniversal}} ($(right 6.north)+(-1,0)$);

        \draw [-Latex] (right B.west) -- ($(right 3.center)+(-1,-2)$);
        \draw [-Latex] (left 3.south-|left A.south) -- node[left]{\cite{Bulatov}, \cite{Zhuk}} (left A); 
        \draw [-Latex] (left 3.north-|left B.north) -- (left B);
        \draw [-Latex] (left B.north) -- (left 2.south-|left B.north);
        \draw [-Latex] (right 2.south-|right A.north) -- (right A);
    \end{tikzpicture}

    \caption{The structure of the Proof. Dashed arrows denote opposite statements. Normal arrows denote implications that will be proved or are already proved. The arrows become a big cycle through all non-dashed boxes when replacing the right side (or left side) by its negations.}
    \label{FigureProofSketch}
\end{figure}
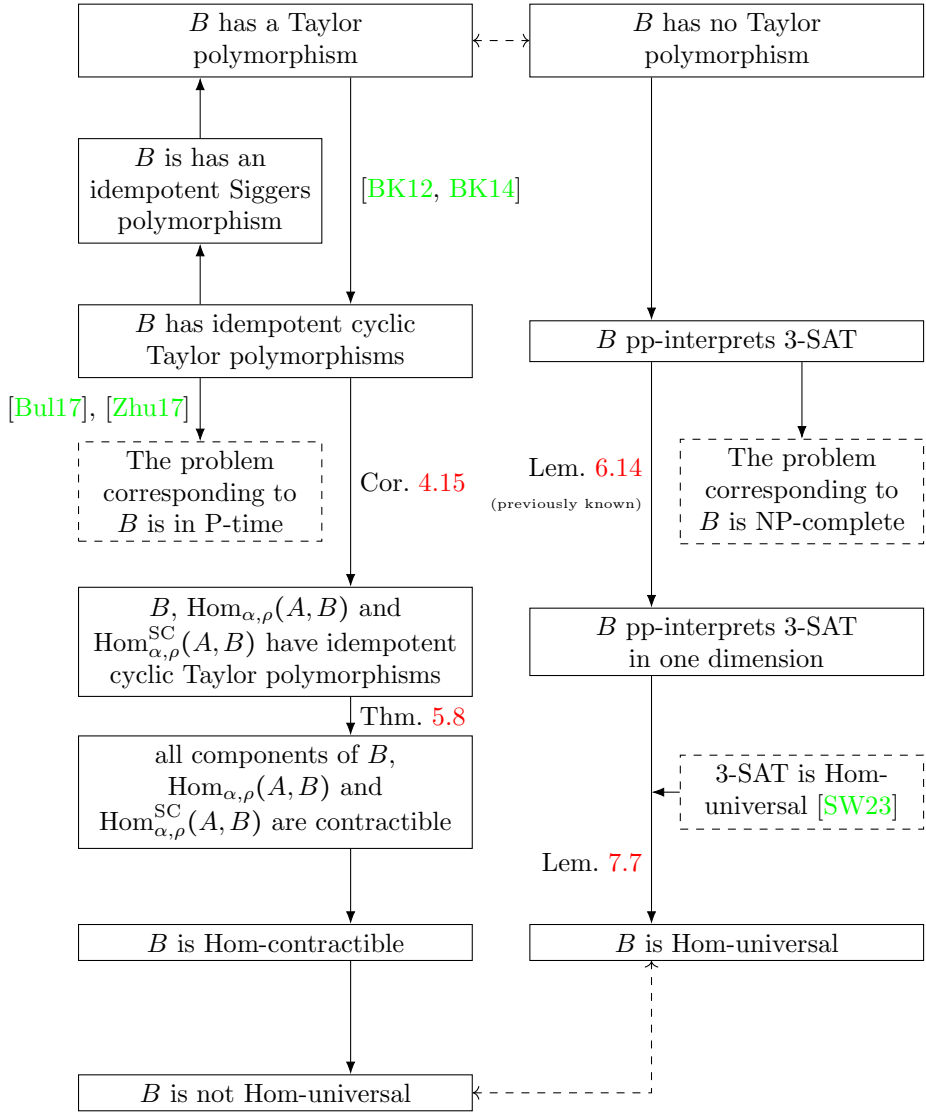

The full structure of the proof is displayed in Figure~\ref{FigureProofSketch}.




\section{Polymorphisms}\label{sectionPolymorphisms}

Polymorphisms are a notion from universal algebra and a generalization of Endomorphisms in categories with products. They give additional characterizations of the dichotomy in Theorem~\ref{TheoremMain}.


\subsection{General Clone Theory}
This section introduces the structure of a clone on the Polymorphisms of an abstract simplicial complex (see Definitions~\ref{DefinitionPolymorphismenmenge} and~\ref{DefinitionCompositionImClone}) and develops some general clone theory. 

\begin{defi} \label{DefinitionPolymorphismenmenge}
    Let $\complexB$ be an abstract simplicial complex. The \de{polymorphisms} of $\complexB$ are given by the set
    $$
        \Pol(\complexB)= \{\clonefunctiona \colon \complexB^n \to \complexB \mid n\in \IN, f \text{ is a simplicial map}\}.
    $$
    A polymorphism $\clonefunctiona$ is \de{idempotent}, if $\clonefunctiona(x,x,\dots,x)=x$ for all $x\in \complexB$. The set of idempotent polymorphisms will be denoted $\Polidem(\complexB)$.
\end{defi}
A similar definition works in other categories with products.
The set of Polymorphisms and the set of idempotent polymorphisms can be naturally seen as clones (or equivalently, as Cartesian operads), where a clone is the following structure:
\begin{defi}
    A \de{clone} is a set $\clone$ which is the disjoint union of a sequence of sets $(\clone_n)_{n\in \IN}$. Moreover, we have for each $n$ exactly $n$ elements $\pi_k^n (1\le k\le n)$ in $\clone_n$ called projections and we have for all $m,n$ composition maps $\circ\colon \clone_n \times (\clone_m)^n\to C_m$ such that
    \begin{align*}
        \forall \clonefunctiona\in \clone_n : &&
        \clonefunctiona \circ (\pi_1^n,...,\pi_n^n) &= \clonefunctiona
        \\
        \forall \clonefunctiona_1,\dots, \clonefunctiona_n\in \clone_n :&&
        \pi_k^n \circ (\clonefunctiona_1,\dots, \clonefunctiona_n) 
        &= \clonefunctiona_k
    \end{align*}
    (projection laws) and
    \begin{align*}
        \forall \clonefunctiona \in \clone_n; \clonefunctionb_1,\dots,\clonefunctionb_n\in \clone_m; \clonefunctionc_1,\dots \clonefunctionc_m \in \clone_k:
        (\clonefunctiona \circ (\clonefunctionb_1,\dots,\clonefunctionb_n)) \circ (\clonefunctionc_1,\dots \clonefunctionc_m)
        \\=
        \clonefunctiona \circ (\clonefunctionb_1\circ (\clonefunctionc_1,\dots,\clonefunctionc_m),\dots,\clonefunctionb_n\circ (\clonefunctionc_1,\dots,\clonefunctionc_m))
    \end{align*}
    (associativity law) hold.
\end{defi}
The representation of a polymorphism as a clone is done in the following way:
\begin{defi}[Polymorphism clone of a simplicial complex] \label{DefinitionCompositionImClone}
    Let $\complexB$ be an abstract simplicial complex. Then, define on $\Pol(\complexB)$ the structure of a clone by defining $\Pol(\complexB)_n=\{\clonefunctiona\colon \complexB^n \to \complexB\}$, $\pi_k^n$ is the projection $\complexB^n \to \complexB$ on the $k$-th component and the composition is given by
    $$
        (\clonefunctiona\circ (\clonefunctionb_1,\dots,\clonefunctionb_n))(x_1,\dots,x_m)=\clonefunctiona(\clonefunctionb(x_1,\dots,x_m),\dots,\clonefunctionb(x_1,\dots x_m)).
    $$
    Define on $\Polidem(\complexB)\subset \Pol(\complexB)$ the substructure.
\end{defi}
It is easy to see that this defines two Clones. 
\begin{defi}
    Let $\clone$ and $\cloneD$ be two clones. Then, a \de{clone homomorphism} from $\clone$ to $\cloneD$ is a map such that $\clone_n$ is mapped inside $\cloneD_n$ and projections $\pi_k^n$ and compositions are preserved.
\end{defi}
There are more important Examples:
\begin{example}
    The \de{projection clone} $\Proj$ is given by $\Proj_n=\{\pi_k^n\mid 1\le k\le n\}$ so it consists only of projections (and all projections are different). The composition follows by the  second projection rule.

    The \de{trivial clone} is a clone $\clone$ where $\clone_n$ consists of only one element.

    The \de{full operation clone on a set} $X$ is a clone $\clone$ where $\clone_n$ consists of all maps $X^n\to X$, $\pi_k^n$ are the projections and $\circ$ is the composition as defined in Definition~\ref{DefinitionCompositionImClone}.
\end{example}
Clearly, the projection Clone maps into every clone and every clone maps into the trivial clone. 
Thus, a clone maps into every other clone if and only if it maps into the projection clone. The distinguishing line if a clone has or has not this property will be important later. We are able to give some criteria already now. As we will discuss later, they turn out to be actually equivalent.

\begin{defi}
    Let $\clone$ be a clone.

    An Element $\clonefunctions$ with the property 
    $$\clonefunctions \circ (\pi_1^3,\pi_1^3,\pi_2^3,\pi_2^3,\pi_3^3,\pi_3^3 ) = \clonefunctions \circ (\pi_3^3,\pi_2^3,\pi_1^3,\pi_3^3,\pi_2^3,\pi_1^3 )$$
    (respectively $\clonefunctions(x,x,y,y,z,z)=\clonefunctions(z,y,x,z,y,x)$ holds for all $x,y,z$)
    is called \de{Siggers element}.

    An element $\clonefunctiona \in \clone_n$ for some $n>1$ such that 
    $$\clonefunctiona = \clonefunctiona \circ (\pi_n^n,\pi_1^n,\pi_2^n,\pi_3^n,\dots,\pi_{n-1}^n) $$
    (respectively $\clonefunctiona(x_1,x_2,\dots,x_n)=\clonefunctiona(x_n,x_1,\dots,x_{n-1})$ for all $x_1,\dots,x_n$) is called \de{cyclic element} of arity $n$.

    More general, an element $\clonefunctionb \in \clone_n$ for some $n>1$ is called \de{Taylor element}, if for every $1\le m\le n$, there is a valid identity of the form
    $$
        \forall x,y : \clonefunctionb(a_1,\dots,a_n)=\clonefunctionb(b_1,\dots,b_n)
    $$
    where $a_i,b_i\in \{x,y\}$ and $a_m=x$ and $b_m=y$. If $\clone$ contains a Taylor element, it is called a \de{Taylor clone}.
\end{defi}
The following remark and lemma are straight forward and general known:
\begin{rem}
    Every clone homomorphism preserves Taylor, Siggers and cyclic polymorphisms.
\end{rem}

\begin{lemma} \label{LemmaEasyTaylor}
    Let $\clone$ be a clone. Then, every item in this list, implies the next one:
    \begin{enumerate}
        \item There is a cyclic element in $\clone$.
        \item There is a Siggers element in $\clone$.
        \item There is a Taylor element in $\clone$.
        \item There is no clone homomorphism $\clone\to \Proj$.
    \end{enumerate}
\end{lemma}


\begin{proof}
We prove the first implication first:

    First note that if $n$ is even, then there is also a cyclic element of arity 2 given by
    $$c \circ  (\underbrace{\underbrace{\pi_1^2,\pi_2^2},\underbrace{\pi_1^2,\pi_2^2}\dots\underbrace{\pi_1^2,\pi_2^2}}_{n/2 \text{ copies of length }2} )$$
    so we may assume $n=2$ or $n$ odd.
    
    For $n=2$, choose $\clonefunctions \coloneqq \clonefunctiona \circ (\pi_2^6, \pi_3^6)$ which equals $\clonefunctiona \circ (\pi_2^2, \pi_1^2) \circ (\pi_2^6, \pi_3^6)=\clonefunctiona \circ (\pi_3^6, \pi_2^6)$ since $\clonefunctiona$ is cyclic. Thus,
    \begin{align*}
        \clonefunctions \circ (\pi_1^3,\pi_1^3,\pi_2^3,\pi_2^3,\pi_3^3,\pi_3^3 ) 
        &=\clonefunctiona \circ (\pi_2^6, \pi_3^6)\circ (\pi_1^3,\pi_1^3,\pi_2^3,\pi_2^3,\pi_3^3,\pi_3^3 )
        \\
        &= \clonefunctiona \circ (\pi_1^3, \pi_2^3)
        \\
        &= \clonefunctiona \circ (\pi_3^6, \pi_2^6) \circ (\pi_3^3,\pi_2^3,\pi_1^3,\pi_3^3,\pi_2^3,\pi_1^3 )
        \\
        &= \clonefunctions \circ (\pi_3^3,\pi_2^3,\pi_1^3,\pi_3^3,\pi_2^3,\pi_1^3 )
    \end{align*}
    For $n$ odd, write $n=2k+3$ (where $k$ can be zero) and define $\clonefunctions$ as
    \begin{align*}
        \clonefunctions 
        &\coloneqq \clonefunctiona \circ (\pi_1^6,\underbrace{\pi_2^6,\dots,\pi_2^6}_{k}, \pi_5^6, \underbrace{\pi_3^6,\dots,\pi_3^6,\pi_3^6}_{k+1})
        \\
        &= \clonefunctiona \circ (\underbrace{\pi_3^6,\pi_3^6,\dots,\pi_3^6}_{k+1}, \pi_1^6,\underbrace{\pi_2^6,\dots,\pi_2^6}_{k}, \pi_5^6)
    \end{align*}
    where the second equality follows because $\clonefunctiona$ is cyclic. 
    In this case, we get
    \begin{align*}
        \clonefunctions \circ (\pi_1^3,\pi_1^3,\pi_2^3,\pi_2^3,\pi_3^3,\pi_3^3 ) 
        &=\clonefunctiona \circ (\underbrace{\pi_1^3,\pi_1^3,\dots,\pi_1^3}_{k+1}, \pi_3^3, \underbrace{\pi_2^3,\dots,\pi_2^3,\pi_2^3}_{k+1})
        \\&= \clonefunctions \circ (\pi_3^3,\pi_2^3,\pi_1^3,\pi_3^3,\pi_2^3,\pi_1^3 )
    \end{align*}
    where we used the definition for the first equal sign and the identity above for the second one.

The second implication is immediate as the Siggers-Identity implies Taylor identities by setting $x=y$, $y=z$ or $z=x$.

For the last implication assume that $\clone$ is a clone with a Taylor element. If there would be a clone Homomorphism $\clone\to \Proj$, then its image would also be a Taylor element. But there are none in $\Proj$, a contradiction.
\end{proof}

There are also a converse result by multiple authors with tame additional assumptions.
\begin{defi}
    An \de{operation clone} is a subset of the full operation clone on a set $X$. An operation clone is called \de{idempotent} if 
    $$
        \forall  \clonefunctiona \in \clone : \clonefunctiona\circ (\pi_1^1,\dots , \pi_1^1)=\pi_1^1
    $$
    holds respectively $\clonefunctiona(x,\dots , x)=x$ holds for all $\clonefunctiona$ and $x$.

    A clone is \de{finitely generated}, if there is a finite subset of $\clone$ such that every element on $\clone$ can be written as a composition of those elements and projections.
\end{defi}

\begin{thm}[Cyclic terms theorem, by \cite{Taylor77} and \cite{CyclicTaylor}] \label{TheoremCyclicTerms}
    Let $\clone$ be an idempotent operation clone on a finite set $X$. Then, the following is equivalent:
    \begin{enumerate}
        \item \label{NrTheoremCyclicTerms1}
        For every prime number $p>|X|$, there is a cyclic element in $\clone_p$.
        \item \label{NrTheoremCyclicTerms2}
        There is a Siggers element in $\clone$.
        \item \label{NrTheoremCyclicTerms3}
        There is a Taylor element in $\clone$.
        \item \label{NrTheoremCyclicTerms4}
        There is no clone homomorphism $\clone\to \Proj$.
    \end{enumerate}    
\end{thm}
\begin{proof}
    The implications $\ref{NrTheoremCyclicTerms1} \implies \ref{NrTheoremCyclicTerms2} \implies \ref{NrTheoremCyclicTerms3} \implies \ref{NrTheoremCyclicTerms4}$ are shown in Lemma~\ref{LemmaEasyTaylor}.
    The backwards implications are more complicated and we give references to the literature.
    
    \begin{description}
    	\item[\ref{NrTheoremCyclicTerms3} $\impliedby$ \ref{NrTheoremCyclicTerms4}] 
    	This is a combination of Corollary 5.2 and Corollary 5.3 from \cite{Taylor77}.
    	\item [\ref{NrTheoremCyclicTerms2} $\impliedby$ \ref{NrTheoremCyclicTerms3}] This was first proven in \cite[Theorem 1.1]{Siggers}.
    	\item [\ref{NrTheoremCyclicTerms1} $\impliedby$ \ref{NrTheoremCyclicTerms3}] This is shown in Theorem 4.1 of \cite{CyclicTaylor} in another language. \qedhere
    \end{description}
\end{proof}
\begin{rem}
	In  \cite{CyclicTaylor} and its reprint \cite{CyclicTaylorReprint}, it is assumed that the clone is finitely generated. This assumption arises from a common notation of a variety and is not needed in the proof, so we can use the more general statement here.
	One can also see 
	by a compactness argument that a clone maps to into $\proj$ if and only if every finitely generated subclone has such a map. Thus, one can derive the more general statement from the restricted version.
\end{rem}

\subsection{Preservation of Homomorphisms and Polymorphisms}

This section is a recall of some facts about abstract simplicial complexes and the Homomorphism functor which have some impact on the theory of Polymorphisms.

\begin{rem}
    Let $\complexA$ and $\complexB$ be abstract simplicial complexes. Let $\subsetA,\subsetAa\subset \complexA$ and $\partialmap\colon \subsetAa \to \complexB$ and map.

    We have $\Hom_{\complexA,\emptyset}(\complexA,\complexB)=\Hom(\complexA,\complexB)$ and $\Homsc_{\complexA,\emptyset}(\complexA,\complexB)=\Homsc(\complexA,\complexB)$. There are natural maps
    $$
    \begin{tikzcd}
        \Hom_{\complexA,\partialmap}(\complexA,\complexB) \ar[r,two heads] \ar[d,hook]
        &
        \Hom_{\subsetA,\partialmap}(\complexA,\complexB) \ar[d,hook]
        \\
        \Hom_{\complexA,\emptyset}(\complexA,\complexB) \ar[r,two heads]
        &
        \Hom_{\subsetA,\emptyset}(\complexA,\complexB)
    \end{tikzcd}
    \text{and}
    \begin{tikzcd}
        \Homsc_{\complexA,\partialmap}(\complexA,\complexB) \ar[r,two heads] \ar[d,hook]
        &
        \Homsc_{\subsetA,\partialmap}(\complexA,\complexB) \ar[d,hook]
        \\
        \Homsc_{\complexA,\emptyset}(\complexA,\complexB) \ar[r,two heads]
        &
        \Homsc_{\subsetA,\emptyset}(\complexA,\complexB)
    \end{tikzcd}
    $$
    given by the restrictions. The surjective maps do not need to be quotient maps as shown in Example~\ref{Example2} and drawn in Figure~\ref{FigureExample2}.
\end{rem}

Also note that most of these constructions behave well with products. This includes the fact that for a simplicial complex, the diagonal map $C\to C^n$ is a simplicial homomorphism and the next lemma:
\begin{lemma}
    Let $\complexA, \complexB$ and $\complexC$ be abstract simplicial complexes. Let $\subsetA, \subsetAa \subset \complexA$ be subsets and $\partialmap_\complexB\colon \subsetAa \to \complexB$ and $\partialmap_\complexC\colon \subsetAa \to \complexC$ be maps on the subset $\subsetAa$.
    
    Then, the following holds:
    \begin{enumerate}
        \item The map 
        \begin{align*}
            \Hom_{\subsetA, \partialmap_\complexB}(\complexA,\complexB) \times \Hom_{\subsetA, \partialmap_\complexC}(\complexA,\complexC) &\to \Hom_{\subsetA, \partialmap_\complexB\times \partialmap_\complexC}(\complexA,\complexB \times \complexC)
            \\
            (f,g) &\mapsto (x\mapsto (f(x),g(x)))        
        \end{align*}
        is an isomorphism of simplicial complexes.
        \item The map 
        \begin{align*}
            \Homsc_{\subsetA, \partialmap_\complexB}(\complexA,\complexB) \times \Homsc_{\subsetA, \partialmap_\complexC}(\complexA,\complexC) &\to \Homsc_{\subsetA, \partialmap_\complexB\times \partialmap_\complexC}(\complexA,\complexB \times \complexC)
            \\
            (f,g) &\mapsto (x\mapsto (f(x),g(x)))
        \end{align*}
        is an isomorphism of simplicial complexes.
    \end{enumerate}
\end{lemma}
\begin{proof}
    By the universal property of the product, this shows the map is a bijection. By the definition of the vertices in $\Hom$ and $\Homsc$, it turns out to be an isomorphism.
\end{proof}

This has an immediate consequence on the polymorphism clones of abstract simplicial complexes, defined in Definitions~\ref{DefinitionPolymorphismenmenge} and~\ref{DefinitionCompositionImClone}:
\begin{lemma} \label{LemmaClonhomo}
    Let $\complexA,$ and $\complexB$ be abstract simplicial complexes. Let $\subsetA, \subsetAa \subset \complexA$ be subsets and $\partialmap\colon \subsetAa \to \complexB$ be a map on the subset $\subsetAa$.
    Then, there are clone homomorphisms
    \begin{align*}
        \Pol(\complexB) &\to \Pol(\Hom_{\subsetA, \partialmap}(\complexA,\complexB)) \text{ by}\\
        \Hom(\complexB^n,\complexB) &\to \Hom((\Hom_{\subsetA, \partialmap}(\complexA,\complexB))^n, \Hom_{\subsetA, \partialmap}(\complexA,\complexB)) \\
        p&\mapsto ((f_1,\dots,f_n)\mapsto (x\mapsto p(f_1(x),\dots,p_n(x))))
    \end{align*}
    and
    \begin{align*}
        \Pol(\complexB) &\to \Pol(\Homsc_{\subsetA, \partialmap}(\complexA,\complexB)) \text{ by}\\
        \Hom(\complexB^n,\complexB) &\to \Hom((\Homsc_{\subsetA, \partialmap}(\complexA,\complexB))^n, \Homsc_{\subsetA, \partialmap}(\complexA,\complexB)) \\
        p&\mapsto ((f_1,\dots,f_n)\mapsto (x\mapsto p(f_1(x),\dots,p_n(x))))
    \end{align*}
    which restrict to clone homomorphisms $\Polidem(\complexB) \to \Polidem(\Hom_{\subsetA, \partialmap}(\complexA,\complexB))$ and $\Polidem(\complexB) \to \Polidem(\Homsc_{\subsetA, \partialmap}(\complexA,\complexB))$.
\end{lemma}
The proof is a simple check of conditions and will be omitted.

This allows us to deduce an implication step of our dichotomy:
\begin{cor} \label{KorollarClonhomozuHom}
    Let $\complexB$ be an abstract simplicial complex such that $\complexB$ has an idempotent cyclic Taylor polymorphism of arity $p$ for infinitely many $p$.

    Then, we get for all abstract simplicial complexes $\complexA$ with subsets $\subsetA,\subsetAa$ and all maps $\partialmap\colon \subsetAa \to \complexB$ that $\Hom_{\subsetA,\partialmap}(\complexA,\complexB)$ and $\Homsc_{\subsetA,\partialmap}(\complexA,\complexB)$ have idempotent cyclic Taylor polymorphism of arity $p$ for infinitely many $p$.
\end{cor}
\begin{proof}
    By Lemma~\ref{LemmaClonhomo}, there are clone homomorphisms from $\Polidem(\complexB)$ to $\Polidem(\Hom_{\subsetA,\partialmap}(\complexA,\complexB))$ and $\Polidem(\Homsc_{\subsetA,\partialmap}(\complexA,\complexB))$. The lemma follows, because clone homomorphisms preserve cyclic terms.
\end{proof}

\section{Polymorphisms and Geometric Properties}
Within the next section, we consider simplicial complexes and want to derive geometric properties from the existence of certain polymorphisms.

\subsection{Some Previous Results on Continuous Polymorphisms of a Simplicial Complex}


There are two meaningful notions of homomorphisms between simplicial complexes: abstract simplicial maps and continuous maps of the geometric representations. 
Both notions extend to different notions of polymorphisms. In our case, the notion of simplicial polymorphism as defined in the Definitions~\ref{DefinitionPolymorphismenmenge} and~\ref{DefinitionCompositionImClone} will be more useful. However, as every simplicial map has a continuous extension to the geometric realization, we should have a look on the more general and more intensely studied case of continuous polymorphisms first to give proper references.


\begin{thm}[{\cite[Corollary 5.3]{Taylor77}}]
    Every finite connected simplicial complex with a continuous idempotent Taylor polymorphism has an abelian fundamental group.
\end{thm}
\begin{rem}
    As every simplicial polymorphism extends to a continuous one, we get also that every finite simplicial complex with a simplicial idempotent Taylor polymorphism has an abelian fundamental group.
\end{rem}

\begin{example} \label{ExampleCircleGroup}
    The unit circle $\IS^1=\{x\in \IC \mid \|x\|=1\}$ has an idempotent continuous Taylor polymorphism $m\colon (\IS^1)^3\to \IS^1, (x,y,z)\mapsto x\cdot y^{-1}\cdot z$. This is a Taylor polymorphism, because $m(x,y,y)=m(y,y,x)=m(x,x,x)=x$ for all $x,y,z$. It is not contractible.
\end{example}

Similar to Example~\ref{ExampleCircleGroup}, every group has an idempotent Taylor polymorphism. And indeed, it holds that:
\begin{cor}
    The fundamental group of a topological group (that is homeomorphic to a simplicial complex) is abelian in every base point.
\end{cor}
The assumption of being homeomorphic to a simplicial complex can be dropped by a different proof using the Eckmann–Hilton argument.

\begin{thm}[various authors] \label{TheoremWhitehead}
    For a finite simplicial complex $\complexA$, the following is equivalent:
    \begin{enumerate}
        \item Every connected component of $\complexA$ is contractible.
        \item The complex $\complexA$ is homotopy equivalent to a discrete set.
        \item It holds $\pi_d(\complexA)=\{0\}$  for all $d\in \IN, n\ge 1$ and all base points. That is
        for all $d\ge 1$, every continuous map  from the $d$-sphere $\IS^d$ to $\complexA$ is homotopic to a constant map.
        (Whitehead Theorem, \cite[Theorem 1]{Whitehead}).
        \item It has a continuous idempotent majority polymorphism, that is a polymorphism $M\colon A^3 \to A$ with $M(x,x,y)=M(x,y,x)=M(y,x,x)=M(x,x,x)=x$ for all $x,y$ (\cite[Theorem 7.7]{Taylor77}).
        \item It has continuous idempotent fully symmetric polymorphisms of all arities, i.e. satisfies
        $$
            \forall n\ge 1 \ \exists P_n\colon \complexA^n \to \complexA \ \forall \sigma \text{ permutation}: P(x_1,\dots,x_n)=P(x_{\sigma(1)},\dots,x_{\sigma(n)})
        $$
            ({\cite{GeneralizedMeans}} and independently {\cite[Theorem 1.1]{SocialChoice}})
    \end{enumerate}
\end{thm}
\begin{rem}
    The theorems are usually stated for the case that $\complexA$ is connected, but this generalization is an easy corollary.

    In the last equivalence, the assumption on $\complexA$ being finite is needed, see \cite{SocialChoice} for a non-finite counterexample.
\end{rem}

\subsection{Some Previous Results on Simplicial Polymorphisms of a Simplicial Complex}

Polymorphisms of abstract simplicial complexes have been studied as well.
We want to mention \cite{TopologyAndAdjunction} in this context. 

Even tough the term Polymorphisms does not appear in that paper, we also want to mention \cite{BooleanCaseLIPIcs} here. In this paper, it is shown, that the simplicial complex which arises from the solutions of a boolean Schäfer constraint satisfaction problem has trivial homology.

However, it was brought to the authors attention, that there is another much earlier paper on finite partially ordered sets, that shows similar results:

\begin{thm}[{\cite[Theorems 2.3 and 3.2]{PosetsLarose}}] \label{TheoremLarose}
    Let $P$ be a finite topological space, not necessarily a Hausdorff space. 
    If $P$ admits a Taylor term, then $P$ is $n$-connected for all $n \ge 0$, i.e. homotopy equivalent to a discrete set.
\end{thm}

The Theorem~\ref{TheoremTaylorContractible} does not follow immediately from this theorem, as the points of a simplicial complex correspond only to the closed points in a finite topological space.
However, there is a proof using similar ideas which will be presented in a separated paper \cite{OprsalMeyer}.

\subsection{Every Taylor Simplicial Complex is Contractible}
The aim of this subsection is to show that every Taylor simplicial complex is contractible, Theorem~\ref{TheoremTaylorContractible}. We give one proof for this theorem that uses the cyclic terms theorem. 

We start by giving an independent elementary proof of the following theorem:
\begin{thm} \label{TheoremSimplicialCyclic}
    Let $\complexA$ be a finite abstract simplicial complex. Assume that there are infinitely many numbers $n$ such that $\complexA$ has a cyclic idempotent polymorphism of arity $n$. Then, every connection component of $\complexA$ is contractible. 
\end{thm}
\begin{proof}[Overview of the proof]
Before going into the details and precise definitions, we give an overview over the proof:
By whiteheads theorem, it suffices to show that each homotopy group is trivial respectively that each map from a sphere into the complex is contractible. By a simplicial approximation, we may assume that the map from the sphere into the complex can be represented by a simplicial map. 

For the 1-sphere, respectively the circle, the simplicial map is given by a map $f\colon \cyclecomplex_n \to \complexA$. By an easy argumentation, we can enlarge $n$ so that we may assume that we have a cyclic polymorphism $c$ of arity $n$. Now observe that the map $\cyclecomplex_n\to \cyclecomplex_n, x \mapsto x + a$ viewed as simplicial map $\IS^1 \to \IS^1$ is homotopic to the identity as it is a rotation. Hence, the map $x \mapsto f(x + a)$ is homotopic to $f$. Now, we get that $g \colon \cyclecomplex_n \to {\complexA}$ defined by
$$
    g(x) = c(f(x), f(x + 1), \dots, f(x + n - 1))
$$
is homotopic to $c(f, f, \dots, f) = f$. It follows immediately that $g$ is a constant map, because $c$ is cyclic. As conclusion, we get that $f$ is homotopic to the constant map $g$ and thus contractible.

For higher spheres, this argument lifts, because $\IS^{n+1}$ is a suspension of $\IS^n$ and thus inductively constructed from $\IS^1$. We therefore can again consider the sphere to be represented by the product of a circle and a hypercube where some circles are contracted to a single point. As above, we consider rotations of the circle and show that $f$ is homotopic to a map $g$ which is constant on all circles and thus contractible.
\end{proof}
\begin{example}
	Consider a simplicial complex with a cyclic polymorphism $c$ of arity 5 and a circle 0--1--2--3--4--0. Then, the drawing
	$$
	\begin{tikzpicture}
		()
		\node (0) at (-5,5) {0} ;
		\node (1) at (-3,5) {1} ;
		\node (2) at (-1,5) {2} ;
		\node (3) at (1,5) {3};
		\node (4) at (3,5) {4};
		\node (5) at (5,5) {0};
		\draw (0)--(1)--(2)--(3)--(4)--(5);
		\node (c0) at (-5,4) {c(0,0,0,0,0)} ;
		\node (c1) at (-3,4) {c(1,1,1,1,1)} ;
		\node (c2) at (-1,4) {c(2,2,2,2,2)} ;
		\node (c3) at  (1,4) {c(3,3,3,3,3)};
		\node (c4) at  (3,4) {c(4,4,4,4,4)};
		\node (c5) at  (5,4) {c(0,0,0,0,0)};
		\draw (c0)--(c1)--(c2)--(c3)--(c4)--(c5);
		\draw[double] (0)--(c0);\draw[double] (1)--(c1);\draw[double] (2)--(c2);
		\draw[double] (3)--(c3);\draw[double] (4)--(c4);\draw[double] (5)--(c5);
		\node (d0) at (-5,3) {c(0,1,1,1,1)} ;
		\node (d1) at (-3,3) {c(1,2,2,2,2)} ;
		\node (d2) at (-1,3) {c(2,3,3,3,3)} ;
		\node (d3) at  (1,3) {c(3,4,4,4,4)};
		\node (d4) at  (3,3) {c(4,0,0,0,0)};
		\node (d5) at  (5,3) {c(0,1,1,1,1)};
		\draw (d0)--(d1)--(d2)--(d3)--(d4)--(d5);
		\draw (c0)--(d0);\draw (c1)--(d1);\draw (c2)--(d2);\draw (c3)--(d3);\draw (c4)--(d4);\draw (c5)--(d5);
		\draw (c1)--(d0);\draw (c2)--(d1);\draw (c3)--(d2);\draw (c4)--(d3);\draw (c5)--(d4);
		\node (e1) at (-5,2) {c(0,1,2,2,2)} ;
		\node (e2) at (-3,2) {c(1,2,3,3,3)} ;
		\node (e3) at (-1,2) {c(2,3,4,4,4)};
		\node (e4) at  (1,2) {c(3,4,0,0,0)};
		\node (e5) at  (3,2) {c(4,0,1,1,1)};
		\node (e6) at  (5,2) {c(0,1,2,2,2)};
		\draw (e1)--(e2)--(e3)--(e4)--(e5)--(e6);
		\draw (e1)--(d1);\draw (e2)--(d2);\draw (e3)--(d3);\draw (e4)--(d4);\draw (e5)--(d5);
		\draw (e1)--(d0);\draw (e2)--(d1);\draw (e3)--(d2);\draw (e4)--(d3);\draw (e5)--(d4);\draw (e6)--(d5);
		\node (f1) at (-5,1) {c(0,1,2,3,3)} ;
		\node (f2) at (-3,1) {c(1,2,3,4,4)} ;
		\node (f3) at (-1,1) {c(2,3,4,0,0)};
		\node (f4) at  (1,1) {c(3,4,0,1,1)};
		\node (f5) at  (3,1) {c(4,0,1,2,2)} ;
		\node (f6) at  (5,1) {c(0,1,2,3,3)} ;
		\draw (f1)--(f2)--(f3)--(f4)--(f5)--(f6);
		\draw (e1)--(f1);\draw (e2)--(f2);\draw (e3)--(f3);\draw (e4)--(f4);\draw (e5)--(f5);\draw (e6)--(f6);
		\draw (e2)--(f1);\draw (e3)--(f2);\draw (e4)--(f3);\draw (e5)--(f4);\draw (e6)--(f5);
		\node (g2) at (-5,0) {c(0,1,2,3,4)} ;
		\node (g3) at (-3,0) {c(1,2,3,4,0)};
		\node (g4) at (-1,0) {c(2,3,4,0,1)};
		\node (g5) at  (1,0) {c(3,4,0,1,2)};
		\node (g6) at  (3,0) {c(4,0,1,2,3)} ;
		\node (g7) at  (5,0) {c(0,1,2,3,4)} ;
		\draw (g2)--(f2);\draw (g3)--(f3);\draw (g4)--(f4);\draw (g5)--(f5);\draw (g6)--(f6);
		\draw (g2)--(f1);\draw (g3)--(f2);\draw (g4)--(f3);\draw (g5)--(f4);\draw (g6)--(f5);\draw (g7)--(f6);
		\draw[double] (g2)--(g3)--(g4)--(g5)--(g6)--(g7);
		
	\end{tikzpicture}
	$$
	gives a homotopy from the circle to the constant on $c(0,1,2,3,4)$. We will define a slightly different homotopy to catch the general case.
\end{example}

As we will have many intermediate results to state during the described steps, we will first prove them independently before collecting the proof of Theorem~\ref{TheoremSimplicialCyclic} at the end of this section to increase the clarity of the proof.

For the first step, recall Whitehead's theorem from Theorem~\ref{TheoremWhitehead} and the definition of a sphere:
\begin{defi}[Sphere]
    Let $\IS^n$ be the $n$-sphere $\{x\in \IR^{n+1} \mid |x|=1\}$ where $|x|$ is the Euclidean absolute value.
\end{defi}


In the next step, we want to represent the continuous map by a nice simplicial map. Thus, we define a nice family of simplicial complexes which are homotopy equivalent to spheres.
\begin{defi}
    Let $H^d_{n,m}$ be the simplicial set inductively defined as follows:
    
    For $d=1$ define $H^1_{n,m}=\cyclecomplex_n$.

    For $d>1$ define $H^{d}_{n,m}=H^{d-1}_{n,m} \times \pathcomplex_m / \sim$ where ${\sim} \subset (H^{d-1}_{n,m} \times \pathcomplex_m)^2$ is the equivalence relation 
    \begin{align*}
        {\sim} = &\{((x,y),(x,y)) \mid (x,y) \in H^{d-1}_{n,m} \times \pathcomplex_m\}
        \\&\cup \{((x,0),(x',0)) \mid x,x' \in H^{d-1}_{n,m} \}
        \\&\cup \{((x,m),(x',m)) \mid x,x' \in H^{d-1}_{n,m} \}
    \end{align*}
    on the vertex set.

    We call $H^d_{n,m}$ the \de{sphere Hypercube Complex}.
\end{defi}
\begin{rem}
    We may consider $H^d_{n,m}$ as quotient of $\pathcomplex_n \times \pathcomplex_m \times \pathcomplex_m \times \dots \times \pathcomplex_m$ where the first path is glued together to a cycle $\cyclecomplex_n$ and for each increasing of $d$ we get two sets which are nailed together to two points.

    So we sometimes consider $H^d_{n,m}$ as product of paths modulo some glued and some nailed points. Note that if a point in $\pathcomplex_n \times \pathcomplex_m \times \pathcomplex_m \times \dots \times \pathcomplex_m$ is nailed, then also all other points which differ only in the first coordinate are again nailed and define the same point.
    
    Note that a subset of a product of paths is a face if and only if in each coordinate the entries of this subset have distance at most one.
\end{rem}
\begin{defi}
    The geometric subdivision $\subdivision{\complexA}$ of an abstract simplicial complex $(\complexA,\facesA)$ is an abstract simplicial complex, whose vertex set is $\facesA$ and a subset $\face$ of $\facesA$ is considered to be a face, if every two elements of $\face$ are comparable by inclusion. We denote the composition of $n$ geometric subdivisions by $\subdivisionn{\complexA}$.
\end{defi}
Recall that every simplicial complex is homeomorphic to its geometric subdivision and recall the next theorem, which can be found for example in \cite[Theorem 2C.1]{Hatcher}:
\begin{thm}[Simplicial Approximation Theorem]
    Consider two simplicial complexes $\complexA$ and $\complexB$. Assume that there is a continuous map $f\colon \georeal{\complexA} \to \georeal{\complexB}$. Then, there is a positive integer $n$ and a simplicial map $g\colon \subdivisionn{\complexA} \to \complexB$ such that the geometric realization $\georeal{g}\colon \georeal{\complexA} \to \georeal{\complexB}$ is homotopic to $f$. Here, we identify $\georeal{\complexA}$ and $\georeal{\subdivisionn{\complexA}}$.
\end{thm}

Now, we can prove the approximation of the spheres by hypercubes:
\begin{lemma} \label{LemmaApproxMoeglich}
The following holds:
\begin{enumerate}
    \item The simplicial set $H^d_{n,m}$ is homotopy equivalent to $\IS^d$ for $m,n\ge 3$.
    \item There is a simplicial homotopy equivalence $H^d_{2dn,2dm} \to \subdivision{H^d_{n,m}}$.
    \item For a  finite abstract simplicial complex $\complexA$ and $d\ge 1$ and a continuous map $m\colon \IS^d \to \georeal{\complexA}$ there are $m$ and $n$, a homotopy equivalence $h\colon \IS^d \to \georeal{H^d_{n,m}}$ and a simplicial map $f\colon H^d_{n,m} \to \complexA$ such that $m=\georeal{f}\circ h$.
\end{enumerate}
\end{lemma}
\begin{proof}
    For the first item, note that $H^d_{n,m}$ is by definition the quotient of $H^{d-1}_{n,m} \times \pathcomplex_m$ where the top and the bottom copy of $H^d_{n,m}$ both get contracted to a single point. On the side of topology, this is called the suspension of $H^d_{n,m}$.

    Recall that $\IS^d$ is a suspension of $\IS^{d-1}$ that is
    \begin{align*}
    	\IS^d &\iso (\IS^{d-1} \times [-1,1] )/ \sim \\
    	{\sim} &=\{(x,x)\} \cup \{(x,y)\mid x,y\in \IS \times \{-1\}\}\cup \{(x,y)\mid x,y\in  \IS \times \{1\} \}
    \end{align*}
    As $H^1_{m,n}$ is a cycle (for $m,n\ge 3$), we get by induction that $H^d_{n,m}$ is homotopy equivalent to $\IS^d$.

    For the second item, note that there is a particular nice map from $\pathcomplex_{2d}^d$ to $\subdivision{\pathcomplex_{1}^d}$.
    This map maps a point $(x_1,x_2,\dots,x_d)$ where each coordinate is divisible by $d$ (so it is $0$, $d$ or $2d$) to the subset of $\subdivision{\pathcomplex_{1}^d}$ that contains $(y_1,y_2,\dots,y_d)$, if and only if $x_i\in \{d, 2dy_i\}$ for all $i$. 
    
    All remaining points are mapped equal as a close point whose coordinated are divisible by $d$. The precise algorithm is the following:
    Consider the smallest $n\ge 0$ such that at most $n$ coordinates are within $[d-n,d+n]$. Threat this coordinates as $d$ and round the others to $0$ respectively $2d$.
    This map is drawn in Figure~\ref{FigureExampleSubdivision} and~\ref{FigureExampleSubdivision2a} for $n\le 3$. As it preserves faces, it extends to a map from $\pathcomplex_{2dm}^{d}$ to $\subdivision{\pathcomplex_{m}^d}$ and from $H^d_{2nd,2md}$ to $\subdivision{H^d_{n,m}}$.

    The third item now follows immediately from the second item and the simplicial approximation theorem. \qedhere

    \begin{figure} 
        \centering
        \begin{tikzpicture}
            \draw (-1,1) node {$\pathcomplex_2$};
            \draw (-3,0) -- (1,0);
            \draw 
                (-3,0) node [fill=white]{$0$}  
                (-1,0) node [fill=white]{$1$}  
                (1,0) node [fill=white]{$2$} ; 
            \draw (5,1) node {$\subdivision{\pathcomplex_1}$};
            \draw (3,0) -- (7,0);
            \draw 
                (3,0) node [fill=white]{$\{0\}$}  
                (5,0) node [fill=white]{$\{0,1\}$}  
                (7,0) node [fill=white]{$\{1\}$} ; 
            \draw (-1,-2) node {$\pathcomplex_4 \times \pathcomplex_4$};
            \draw (-3,-3) -- (1,-3) -- (1,-7) -- (-3,-7) -- cycle;
            \draw 
                (-3,-3) node [fill=white]{$\{00\}$}
                (-3,-4) node [fill=white]{$\{00\}$}
                (-2,-3) node [fill=white]{$\{00\}$}
                (-2,-4) node [fill=white]{$\{00\}$};
            \draw (-3,-3) -- (1,-3) -- (1,-7) -- (-3,-7) -- cycle;
            \draw (-3,-4) -- (1,-4) -- (1,-6) -- (-3,-6) -- cycle;
            \draw (-2,-3) -- (0,-3) -- (0,-7) -- (-2,-7) -- cycle;
            \draw (-1,-3) -- (-1,-7);
            \draw (-3,-5) -- (1,-5);
            \draw
                (-0,-3) -- (1,-4)
                (-1,-3) -- (1,-5)
                (-2,-3) -- (1,-6)
                (-3,-3) -- (1,-7)
                (-3,-4) -- (0,-7)
                (-3,-5) -- (-1,-7)
                (-3,-6) -- (-2,-7)
                ;
            \draw 
                (-2,-3) -- (-3,-4)
                (-1,-3) -- (-3,-5)
                (0,-3) -- (-3,-6)
                (1,-3) -- (-3,-7)
                (1,-4) -- (-2,-7)
                (1,-5) -- (-1,-7)
                (1,-6) -- (0,-7)
                ;
            \draw 
                (-3,-3) node [fill=white]{$\{00\}$}
                (-3,-4) node [fill=white]{$\{00\}$}
                (-2,-3) node [fill=white]{$\{00\}$}
                (-2,-4) node [fill=white]{$\{00\}$}
                (0,-3) node [fill=white]{$\{01\}$}
                (0,-4) node [fill=white]{$\{01\}$}
                (1,-3) node [fill=white]{$\{01\}$}
                (1,-4) node [fill=white]{$\{01\}$}
                (-3,-6) node [fill=white]{$\{10\}$}
                (-3,-7) node [fill=white]{$\{10\}$}
                (-2,-6) node [fill=white]{$\{10\}$}
                (-2,-7) node [fill=white]{$\{10\}$}
                (0,-6) node [fill=white]{$\{11\}$}
                (0,-7) node [fill=white]{$\{11\}$}
                (1,-6) node [fill=white]{$\{11\}$}
                (1,-7) node [fill=white]{$\{11\}$}
                (-1,-3) node [fill=white]{$\{00,01\}$}
                (-1,-4) node [fill=white]{$c$}
                (-1,-6) node [fill=white]{$c$}
                (-1,-7) node [fill=white]{$\{10,11\}$}  ;
            \draw 
                (-3,-5) node [align=center, fill=white]{$\{00,$\\$10\}$}
                (-2,-5) node [align=center, fill=white]{$c$}
                (-1,-5) node [align=center, fill=white]{$c$}
                (-0,-5) node [align=center, fill=white]{$c$}
                (+1,-5) node [align=center, fill=white]{$\{01,$\\$11\}$};
            \draw (+5,-2) node [align=center]{A subcomplex of \\ $\subdivision{\pathcomplex_1 \times \pathcomplex_1}$};
            \draw (3,-3) -- (7,-3) -- (7,-7) -- (3,-7) -- cycle;
            \draw
                (3,-3) -- (7,-7)
                (3,-7) -- (7,-3)
                (5,-3) -- (5,-7)
                (3,-5) -- (7,-5);
            \draw
                (3,-3) node [fill=white]{$\{00\}$}
                (7,-3) node [fill=white]{$\{01\}$}
                (5,-3) node [fill=white]{$\{00,01\}$}
                (7,-7) node [fill=white]{$\{11\}$}
                (3,-7) node [fill=white]{$\{10\}$}
                (5,-7) node [fill=white]{$\{10,11\}$}  
                (3,-5) node [align=center, fill=white]{$\{00,$\\$10\}$}
                (7,-5) node [align=center, fill=white]{$\{01,$\\$11\}$}
                (5,-5) node [align=center, fill=white]{$c\coloneqq $\\$\{00,01,$\\$10,11\}$}
            ;
        \end{tikzpicture}
        \caption{The path $\pathcomplex_2$ is isomorphic to $\subdivision{\pathcomplex_1}$ and $\pathcomplex_4 \times \pathcomplex_4$ maps into $\subdivision{\pathcomplex_1 \times \pathcomplex_1}$. To denote this map, the points of $\pathcomplex_4 \times \pathcomplex_4$ are labeled with their respective images in $\subdivision{\pathcomplex_1 \times \pathcomplex_1}$.} \label{FigureExampleSubdivision}
    \end{figure}
    \begin{figure}
        \centering \tabcolsep=2px
        \begin{tabular}{ccccccc|ccccccc|ccccccc|ccccccc}
             (1)&1 & 1 &(2)& 1 & 1 &(1)  & 1 & 1 & 1 & 2 & 1 & 1 & 1   & 1 & 1 & 1 & 4 & 1 & 1 & 1   &(2)& 2 & 4 &(4)& 4 & 2 &(2)\\
             1 & 1 & 1 & 2 & 1 & 1 & 1   & 1 & 1 & 1 & 2 & 1 & 1 & 1   & 1 & 1 & 1 & 8 & 1 & 1 & 1   & 2 & 2 & 8 & 8 & 8 & 2 & 2 \\
             1 & 1 & 1 & 4 & 1 & 1 & 1   & 1 & 1 & 1 & 8 & 1 & 1 & 1   & 1 & 1 & 1 & 8 & 1 & 1 & 1   & 4 & 8 & 8 & 8 & 8 & 8 & 4 \\
             (2)&2 & 4 &(4)& 4 & 2 &(2)  & 2 & 2 & 8 & 8 & 8 & 1 & 1   & 4 & 8 & 8 & 8 & 8 & 8 & 4   &(4)& 8 & 8 &(8)& 8 & 2 &(4)\\
             1 & 1 & 1 & 4 & 1 & 1 & 1   & 1 & 1 & 1 & 8 & 1 & 1 & 1   & 1 & 1 & 1 & 8 & 1 & 1 & 1   & 4 & 8 & 8 & 8 & 8 & 8 & 4 \\
             1 & 1 & 1 & 2 & 1 & 1 & 1   & 1 & 1 & 1 & 2 & 1 & 1 & 1   & 1 & 1 & 1 & 8 & 1 & 1 & 1   & 2 & 2 & 8 & 8 & 8 & 2 & 2 \\
             (1)&1 & 1 &(2)& 1 & 1 &(1)  & 1 & 1 & 1 & 2 & 1 & 1 & 1   & 1 & 1 & 1 & 4 & 1 & 1 & 1   &(2)& 2 & 4 &(4)& 4 & 2 &(2)\\
             \hline
             \multicolumn{7}{c|}{layer 0 and 6} & \multicolumn{7}{c|}{layer 1 and 5} & \multicolumn{7}{c|}{layer 2 and 4} & \multicolumn{7}{c}{layer 3} 
        \end{tabular}
        \caption{The map for $d=3$ is a map from $\pathcomplex_{6}^3$ to $\subdivision{\pathcomplex_1^3}$. 
            To safe some paper and clarity, 
            we only display the number of elements, the image has. So the number 8 which is in brackets represents the fact that the vertex $(3,3,3)$, which is the center of the (hyper-)cube, is mapped to the only subset of $\{0,1\}^3$ with 8 elements. The number at a position $(x_1,x_2,x_3)$ is displayed in brackets, if all coordinates $x_1$ to $x_3$ are $0$, $3$ or $6$.}
        \label{FigureExampleSubdivision2a}
    \end{figure}
\end{proof}

Now, we have in fact reduced Theorem~\ref{TheoremSimplicialCyclic} to showing that every simplicial map $H^d_{m,n} \to \complexA$ is contractible. This will be shown in the next lemma:

\begin{lemma} \label{LemmaMapContractible}
    Let $\complexA$ be an abstract simplicial complex and assume that $\complexA$ has a $p$ dimensional cyclic Taylor polymorphism for infinitely many $p$.
    Let $H^d_{n,m}$ be a sphere hypercube complex and $f$ be a simplicial map $f\colon H^d_{n,m} \to \complexA$.

    Then, $f$ is homotopic to a constant map $f_1$. 
\end{lemma}

\begin{proof}
    We want to prove this lemma in three steps by enlarging the hypercube, constructing another map $f_2\colon H^d_{n,m} \to \complexA$ such that $f$ is homotopic to $f_2$ where the homotopy fixes all nailed point and then proving that $f_2$ factors through a contractible complex.
    
    To enlarge the sphere hypercube complex. We write $(x,\yvect)$ for a point with first coordinate $x$ and other coordinates $\yvect$. Note that the first dimension of a point is indeed special in the sense that if a point $(x,\yvect)$ is nailed, every point $(x',\yvect)$ is also nailed and they even represent the same point. 
    To enlarge the sphere hypercube complex, we want to add an additional $(d-1)$-dimensional hyperplane which is orthogonal to the first component between two layers $i$ and $i+1$ respectively between $\{(i,\yvect)\mid \yvect\}$ and $\{(i+1,\yvect)\mid \yvect\}$.
    Given $H^d_{n,m}$ and a simplicial map $f\colon H^d_{n,m}\to \complexA$, consider the complex $H^d_{n+1,m}$ and as map $g\colon H^d_{n+1,m}\to \complexA$ we consider any of the maps such that
    $$
    	g(x,\yvect)\in \begin{cases}
    		\{f(x,\yvect)\} &\text{if } x\le i
    		\\
    		\{f(x-1,\yvect), f(x,\yvect)\} &\text{if } x= i+1
    		\\
    		\{f(x-1,\yvect)\} &\text{if } x> i+1
    	\end{cases}
    $$
    holds. Then, a careful look at the definitions show that each of this maps indeed preserves glued and nailed points. They are simplicial and contractible if and only if $f$ is. As we want to repeat this process a cuple of times, we rename $n+1$ to $n$ and $g$ to $f$ after the process.
    
    By repeating this process, we may assume that if the map $f$ changes in the first component at a position, $f(x,\yvect)\ne f(x+1,\yvect)$, then for each other  point $(x,\yvect')$ with the same first coordinate, the values of $f(x,\yvect')$ and $f(x+1,\yvect')$ coincide. This can simply be done by adding enough different hyperplanes between the original planes at $x$ and $x+1$. Note that this property cannot be undone by further enlarging the hypercube.
    
    However, by further enlarging the cube, we may assume that two of three consecutive hyperplanes orthogonal to the first dimension are actually equal. We call $i$ such that the hyperplanes $i$ and $i-1$ are equal an equallity number. By adding additional hyperplanes between two equal hyperplanes, we may furthermore choose the number $n$ of vertices in the first direction such that $\complexA$ carries a cyclic polymorphism $p$ of arity $n$.

    From now, $n$ will refer to the new value. We have created a sphere hypercube complex $H=H^d_{n,m}$ and a map $f\colon H\to \complexA$ which is homotopic to a constant map, if and only if the $f$ from the beginning was. 
    We are now able to construct the map $f_2$ and the homotopy to it.
    For doing this, we switch into maps $H\to \complexA^n$.
    
    Let $g^i$ with $i\in \IN, 0\le i\le 2n$ denote a couple of maps from $H$ to $\complexA^n$. With $g^i_j$, $1\le j\le n$, we denote the $j$-th component, which is a map $H\to \complexA$.
     
    Define $g^0$ to be the map given by $f$ followed by the diagonal inclusion $\Delta_n\colon \complexA\to \complexA^n, z\mapsto (z,...,z)$. Define more generally $g^{2i}_j$ to be the map
    \begin{align*}
    	g^{2i}_j\colon H &\to \complexA \\
    	(x,\yvect)&\mapsto
    	\begin{cases}
    		f(x+j,\yvect) &\text{if }j\le i\\
    		f(x+i,\yvect) &\text{if }j> i
	    \end{cases}
    \end{align*}
	where $x+i$ is considered modulo $n$.
	Define intermediate functions
	\begin{align*}
		g^{2i+1}_j\colon H &\to \complexA \\
		(x,\yvect)&\mapsto
		\begin{cases}
			g^{2i+2}_j &\text{if $x$ is an equality number}\\
			g^{2i}_j &\text{else}
		\end{cases}
	\end{align*}
	Now, the maps $(g^i)_{0\le i\le 2n}$ define a homotopy from $g^0$ to $g^{2n}$ in the even stronger sense that the map
	\begin{align*}
		G\colon H \times P_{2n} \to \complexA^n \\
		((x,\yvect), i)\mapsto g^i(x,\yvect)
	\end{align*}
	is simplicial. 
 Note furthermore that if $(x,\yvect)$ is a nailed point, then $g^{i}(x,\yvect)$ is constantly $f(x,\yvect)$. 
    
    Recall that $p$ is a simplicial homomorphism $\complexA^n \to \complexA$. Thus we can compose it and look at $p\circ G$. This defines a homotopy from $p\circ g^0=p\circ \Delta_n\circ f=f$ to $p\circ g^{2n}$ as $p$ is idempotent. We name $f_2\coloneqq p\circ g^{2n}$.
    
    Note that 
    $$
    	f_2(x,\yvect)=p(f(x,\yvect),f(x+1,\yvect),...,f(x+n-1,\yvect))
    $$
    is independent from $x$ as $p$ is cyclic. Thus, let $I$ be the quotient of $H$ where all points are identified if they only differ in the first coordinate. Then, the map $f_2$ factors through $I$. 
    
    Note that $I$ is contractible. This is best visible on the topological side. While $H$ represented $\IS^d$ which is $\IS^1$ suspended $d-1$ times, $I$ represents the space given by a point suspended $d-1$ times. This space is contractible by induction as a point is contractible and a suspension of a contractible space is contractible. 
    Thus, $f_2$ factors through a contractible simplex $I$ and thus has a homotopy to a constant map $f_1$.
    
    Therefore, also every $f\colon \IS^d \to \complexA$ is homotopic to a constant map $f_1$.
\end{proof}

\begin{proof}[Full proof of Theorem~\ref{TheoremSimplicialCyclic}]
    Consider any continuous map $m\colon \IS^d \to \complexA$ where $d\ge 1$. By Lemma~\ref{LemmaApproxMoeglich}, there is a simplicial complex $H^d_{n,m}$ a continuous map $g\colon \IS^d\to H^d_{n,m}$ and a simplicial map $f\colon H \to \complexA$ such that $m=fg$. By Lemma~\ref{LemmaMapContractible}, the map $f$ is contractible, so $m$ is also contractible.
\end{proof}

\begin{thm} \label{TheoremTaylorContractible}
    Let $\complexA$ be a finite abstract simplicial complex that has a simplicial idempotent Taylor polymorphism. Then, every connected component of $A$ is contractible.
\end{thm}
\begin{proof}[First proof]
    If $\complexA$ has an idempotent Taylor polymorphism, then by Theorem~\ref{TheoremCyclicTerms}, it has cyclic maps of infinite arities among the clone of idempotent polymorphisms. By Theorem~\ref{TheoremSimplicialCyclic}, every component of $\complexA$ is contractible.
\end{proof}
\begin{proof}[Second proof]
	There is an alternative proof generalizing the ideas from Theorem~\ref{TheoremLarose} which is not relying on the cyclic terms theorem. We will give this proof in a separate paper \cite{OprsalMeyer}.
\end{proof}
Note that Theorem~\ref{TheoremTaylorContractible} together with~\ref{KorollarClonhomozuHom} imply that every Taylor simplicial complex is \contractible{}.
\section{Relational Structures}

In this part of the paper, we want to show that a finite complex $\complexA$ without a Taylor polymorphism or equivalently with a clone homomorphism $\Pol(\complexA)\to \proj$ is \universal{}. To do this, we want to use in a first step the clone homomorphism and Birkhoff's HSP Theorem. We introduce relational structures as they contain abstract simplicial complexes and are a convenient setting for Birkhoff's Theorem.

\subsection{Basic Definitions of Relational Structures}

We start by defining relational structures.
\begin{defi}
    A \de{relational signature} $\signature$ consists of a set (again denoted by $\signature$) called relation symbols and a map called arity from this set to the positive integers.

    A \de{relational structure} $\structA$ over a relational signature $\signature$ or short a \de{$\signature$-structure} is a set $\setstructA$ and for each relation symbol $\Relationsymbol\in \signature$ of arity $n$, there is a subset $\Relationsymbol^\structA \subset \setstructA^n$.

    We call a $\signature$-structure $\structA$ \de{finite}, if both $\signature$ and $\setstructA$ are finite. Unless specified otherwise, we denote the base set of $\structA$ again by $\structA$ instead of $\setstructA$.
\end{defi}
There are also the notion of homomorphisms and polymorphisms of structures:
\begin{defi}
    For two $\signature$-structures $\structA$ and $\structB$, the set 
    $$
        \Hom_\signature(\structA,\structB)=\{f\colon \setstructA \to \setstructB \mid \forall \Relationsymbol \in \sigma : f(\Relationsymbol^\structA)\subset \Relationsymbol^\structB  \}
    $$
    are the $\signature$-\de{homomorphisms}.
    The \de{polymorphism clone} of $\structB$ is the set of the $\signature$-homomorphisms $\structB^n\to \structB$ where the composition is defined as in Definition~\ref{DefinitionCompositionImClone}.
\end{defi}

Many objects from algebra or combinatorics can be considered as relational structures. 

\begin{defi}
    Let $\dSAT$ be the relational structure on the set $\{0,1\}$ over the signature $\{0,1, \NOT, \dOR\}$ where the relations $0$ and $1$ contain their respective Element,
    \begin{align*}
    {\NOT} &= \{(0,1),(1,0)\} \text{ and } \\ {\dOR} &= \{0,1\}^3\setminus \{(0,0,0)\}.
    \end{align*}
\end{defi}
The structure $\dSAT$ will play an important role in this paper.

\begin{lemma}[{\cite[Example 33]{PolyUse}}] \label{LemmaDSATproj}
    We have $\Pol(\dSAT)=\proj$. 
\end{lemma}

\subsection{Simplicial Complexes as Relational Structures}
The simplest way to make an abstract simplicial complex a relational structure is the following: 
\begin{defi}
    Let $\signature$ be the set $\{\facerel_n \mid n\in \IN \cup \{0\}\}$ where $\facerel_n$ is a relation symbol of size $n+1$. 
    Then, an abstract simplicial complex $(\complexB,\facesB)$ induces a $\signature$-structure $\structB$ by 
    setting $\facerel_n$ to all faces of size $n$ respectively by defining $\facerel_n^\structB=\{(x_0,\dots , x_n)\mid \{x_0,\dots, x_n\}\in \facesB\}$.
\end{defi}
The simplicial homomorphisms between two abstract simplicial complexes are exactly the $\signature$-homomorphisms between the induced structures. However, as $\signature$ is always infinite, we will not be able to use results obtained for finite structures. This, however can be solved by a small trick using the next lemma:
\begin{lemma} \label{LemmaDimensionForHomo}
    Let $\complexA$ and $\complexB$ be two simplicial complexes such that the dimension $\dimension$ of $\complexB$ is finite. Let $f\colon \complexA \to \complexB$ be a map of sets such that $f$ preserves the relations $\facerel_n$ for all $n\le \dimension+1$. Then, $f$ is simplicial.
\end{lemma}
\begin{proof}
    Assume that $f$ is no simplicial map. Then, there is a face $\face\subset \complexA$ such that $f(\face)$ is not a face. As every face of $\complexB$ has at most $d+1$ vertices, there is a subset $C\subset f(\face)\subset \complexB$ which is not a face and has size $k+1\le d+2$. 
    Let $\{x_0,\dots,x_k\}\coloneqq D\subset \complexA$ be a set which contains a single preimage for each element in $C$. Then, the set $D$ is a subset of a face and thus itself a face. Therefore, the tuple $(x_0,\dots,x_k)$ is in $\facerel_k^\structA$. 
    That is a contradiction as $f$ preserves $\facerel_k$ but $(f(x_0),\dots, f(x_k))$ cannot be in $\facerel_k^\structB$ as the set $\{f(x_0),\dots ,f(x_k)\}=C$ is not a face.
\end{proof}
\begin{rem}
    Note that in the previous lemma, the bound $n\le d+1$ is a crucial difference in the dimension. We need to consider for example the 1-cells of $\complexA$ if $\complexB$ is a collection of points and thus has dimension 0.
\end{rem}

Now, we can consider an $\dimension$-dimensional complex over the finite signature $\{\facerel_n\mid n\le d+1\}$. As every finite abstract simplicial complex is also finite dimensional, this makes a finite abstract simplicial complex into a finite relational structure. However, there is a slight difference in the definition, we will make. The reason for this adjustment is that complexes have a lot of endomorphisms including constant endomorphisms which we want to avoid. This adjustment corresponds to the partially defined map $\partialmap$ in $\Homsc_{\subsetA,\partialmap}$.

\begin{defi}
    Let $(\complexB,\facesB)$ be a finite simplicial complex with dimension $\dimension$. Define the \de{finite idempotent realization} of $\complexB$ as follows: Take the signature $\{\facerel_n\mid n\in \IZ, 0\le n\le \dimension+1\}\cup \{\Relationsymbol_x \mid x\in \complexB\}$ where the arity of $\facerel_n$ is $n+1$ and the arity of $\Relationsymbol_x$ is 1. Now, the finite realization is the structure $\structB$ over this signature where the set is given by $\complexB$ and the relations by
    $$(x_0,\dots ,x_n)\in \facerel_n^\complexB \iff \{x_0,\dots , x_n\}\in \facesB $$
    and $\forall x \in \complexB : \Relationsymbol_x^\complexB = \{x\}$.
    
    We call a relational structure \de{idempotent} if the identity is the only endomorphism.
\end{defi}

\begin{lemma} \label{LemmaPolErhalten}
    Let $(\complexB,\facesB)$ be a finite simplicial complex and let $\structB$ be its finite idempotent realization. Then, $\Polidem(\complexB) = \Pol(\structB)$.
\end{lemma}
\begin{proof}
	A map $\complexB^n \to \complexB$ is simplicial if and only if it preserves the relations $\{\facerel_n\mid n\in \IZ, 0\le n\le \dimension+1\}$ by Lemma~\ref{LemmaDimensionForHomo}. It moreover is idempotent, if and only if it preserves every point that is if it preserves  the relations $\{\Relationsymbol_x \mid x\in \complexB\}$.
\end{proof}

\subsection{Primitive Positive Interpretations}
Primitive positive interpretations are closely connected to clone homomorphisms and make the existence of a clone homomorphism applicable.

\begin{defi}
    Let $\signature$ be a relational signature. A \de{primitive positive $\signature$-formula} is a first order formula which is recursively build by
    \begin{itemize}
        \item $\exists x: \phi$ where $\phi$ is again a primitive positive $\signature$-formula
        \item $\phi \land  \psi$ where $\phi$ and $\psi$ are again primitive positive $\signature$-formulas
        \item $x=y$ where $x$ and $y$ are variables
        \item $R(x_1,\dots,x_n)$ where $R$ is an $n$-ary relation from $\signature$ and $x_1$ to $x_n$ are variables
        \item true, the formula which is always valid and
        \item false, the formula which is never valid
    \end{itemize}
    A primitive positive formula can have free variables.
    
    Let $\structA$ be a $\signature$-structure and $\structB$ be a $\tau$ structure. Then, $\structA$ \de{primitively positively interprets} $\structB$ 
    if and only if there is a positive integer $d$ and there are primitive positive formulas $\psi$ with $d$ free variables, $\phi_=$ with $2d$ free variables and for each relation symbol $R$ in $\tau$ of arity $n$ a primitive positive formula $\phi_R$ with $nd$ free variables and a surjective map
    $$
        h\colon \{(x_1,\dots x_d)\in \structA^d \mid \psi(x_1,\dots x_d)\} \to \structB
    $$
    such that the equivalence
    $$
        h(x_1,\dots x_d) = h(x'_1,\dots x'_d) \iff \phi_=(x_1,\dots x_d,x'_1,\dots x'_d)
    $$
    and for all relation symbols $R\in \tau$ of arity $n$ the equivalence
    $$
        R(h(x_{1,1},\dots x_{1,d}), \dots,h(x_{n,1},\dots, x_{n,d} ))
        \iff
        \psi_{R}(x_{1,1},\dots x_{1,d}, \dots,x_{n,1},\dots, x_{n,d} )
    $$
    hold. The number $d$ is called \de{dimension} of the interpretation.

    We commonly use \de{pp-} instead of primitive positive.
\end{defi}

\begin{thm}[Inv-Pol-Theorem, \cite{GeigerInvPol} and \cite{BKKRInvPol}] 
	\label{TheoreminvPolA}
    Let $\structB$ be a $\signature$-structure and $\subsetB \subset \setstructB^d$ be a subset of a power. Then, the following are equivalent:
    \begin{enumerate}
        \item There is a primitive positive $\signature$-formula $\phi$ with $d$ many free variables such that for all $x_1,\dots,x_d$ in $\setstructB$, we get $\phi(x_1,\dots,x_d) \iff (x_1,\dots,x_d)\in \subsetB$.
        \item The set $\subsetB$ is invariant under polymorphisms, that is for all $f\in \Pol(\structB)$ of any degree $n$, and all $x^1,\dots,x^n \in \subsetB\subset \setstructB^d$ collection of $n$ many $d$-tuples, the image
        $(f(x^1_1,\dots x^n_1),\dots,f(x^1_d,\dots x^n_d))$ is again in $\subsetB$.
    \end{enumerate}
\end{thm}

\begin{rem} \label{RemarkPhiA}
    Let $\structA$ be a finite idempotent relational structure and $a\in \structA$. Then, there is by Theorem~\ref{TheoreminvPolA} a pp-formula $\phi_a$ with a single free variable such that $\{a\}=\{x\in \structA \mid \phi_a(x)\}$.
\end{rem}
Some authors also consider the existence of $\phi_a$ or even that $\phi_a$ is given by a single unary relation as part of the definition of an idempotent structure.

We give Birkhoff's HSP theorem in terms of structures.
\begin{thm}[Birkhoff's HSP theorem \cite{BirkhoffHSP}]
	
	\label{TheoremPolPPInterpret}
	Let $\structA$ be a finite $\signature$-structure and $\structB$ be a finite $\tau$ structure.
	
	Then, $\structA$ pp-interprets $\structB$ if and only if there is a clone Homomorphism $\Pol(\structA) \to \Pol(\structB)$.
\end{thm}

The next lemma is already known in the constraint satisfaction problem  community. However, since it is an important step to apply Lemma~\ref{LemmaPPIntHomotopie} later, we include a proof.
\begin{lemma}\label{LemmaConsturctInterpret}
    Let $\structA$ be a finite structure. The following are equivalent:
    \begin{enumerate}
        \item \label{NrLemmaConstructionInterpret1}
        There is a clone homomorphism $\Pol(\structA)\to \proj$.
        \item \label{NrLemmaConstructionInterpret2}
        The structure $\structA$ pp-interprets 3-SAT.
        \item \label{NrLemmaConstructionInterpret3}
        The structure $\structA$ pp-interprets 3-SAT with a 1-dimensional interpretation.
    \end{enumerate}
\end{lemma}
\begin{proof}
   \ref{NrLemmaConstructionInterpret1} $\implies$ \ref{NrLemmaConstructionInterpret2}:
    If there is a clone Homomorphism from $\Pol(\structA)$ into $\Proj$ then there is a clone homomorphism from $\Pol(\structA)$ into every clone, including $\Pol(\dSAT)$. By Theorem~\ref{TheoremPolPPInterpret}, $\structA$ pp-interprets $\dSAT$. 

   \ref{NrLemmaConstructionInterpret2} $\implies$ \ref{NrLemmaConstructionInterpret3}:
    Choose one pp-interpretation.
    Let $d$, $\psi$, $\phi_=$, $\phi_{\NOT}$, $\phi_{\dOR}$ and $h$ be the objects from this interpretation. 
    We look at the $h$-preimages of $\{0\}$ and $\{1\}$ in $\structA^d$ and focus on the first coordinate. Note that $d\ge 1$, because $\structA^0$ has only one element.
    
    If there is any $a \in \structA$ such that there exist $x_2,\dots,x_n$ and $x'_2,\dots,x'_n$ in $\structA$ such that $h(a,x_2,\dots,x_n)=0$ and $h(a,x'_2,\dots,x'_n)=1$ then we can replace $\psi(x_1,\dots,x_n)$ by $\psi(x_1,\dots,x_n) \land \phi_a(x_1)$ where $\phi_a$ is as in Remark~\ref{RemarkPhiA}. This results in a second $d$-dimensional pp-interpretation.
    However, this second pp-interpretation behaves in fact like a $d-1$ dimensional interpretation because $x_1$ will always take the same value. Therefore, we can define a $d-1$ dimensional interpretation by replacing $x_1$ with a new existentially quantified variable. This decreases $d$ by one.

    In the other case, there is no $a \in \structA$ such that $(a,x_2,\dots,x_n)$ can be mapped to both $0$ and $1$. In this case, we can existentially quantify every other variable as the first one already contains all the information. This defines a second pp-interpretation with $d=1$.

    We can inductively repeat this construction until $d=1$. Therefore, there is always a pp-interpretation in this dimension.

   \ref{NrLemmaConstructionInterpret3} $\implies$ \ref{NrLemmaConstructionInterpret2}: trivial.

   \ref{NrLemmaConstructionInterpret2} $\implies$ \ref{NrLemmaConstructionInterpret1}: This follows from $\Pol(\dSAT)\iso \proj$ and Theorem~\ref{TheoremPolPPInterpret}.
\end{proof}

    
    
    
    

\subsection{Connections to Complexity Theory}
To each relational structure, there exists a related decision problem where the complexity of the decision problem depends on algebraic conditions of the structure.
\begin{defi}
	The \de{Constraint satisfaction Problem (CSP)} of a structure $\structB$ is to decide for any finite structure $\structA$ over the same signature whether there is a homomorphism $\structA \to \structB$.
\end{defi}
\begin{thm}[Bulatov \cite{Bulatov} and independently Zhuk \cite{Zhuk}]\label{TheoremZhuk}
	Consider a finite structure $\structB$.
	If $\structB$ has a Taylor Polymorphisms, then $\CSP(\structB)$ is in P. Otherwise, $\CSP(\structB)$ is NP-complete.
\end{thm}
We can adopt this setting to simplicial complexes.
\begin{defi}
	For an abstract simplicial complex $\complexB$, its \de{precolored CSP} is the problem which gets as input another simplicial complex $\complexA$ (given by the number of vertices followed by the list of all faces), a subset $\subsetAa \subset \complexA$ and a map $\partialmap\colon\subsetAa \to \complexB$ and has to decide if there is a homomorphism $\complexA \to \complexB$ that agrees with $\partialmap$ on $\subsetAa$ i.e. if $\Hom_{\complexA, \partialmap}(\complexA,\complexB)$ is nonempty.
\end{defi}

\begin{lemma} \label{LemmaPReduktion}
	Let $(\complexB,\facesB)$ be a finite simplicial complex and let $\structB$ be its finite idempotent realization. Then, the precolored CSP of $\complexB$ and the CSP of $\structB$ are polynomial time equivalent.
\end{lemma}
\begin{proof}
	We prove both directions separately.
	
	Given any relational structure $\structA=(\setstructA; \facerel_n^\structA \mid n\le \dim(\complexB)+1, \Relationsymbol_x^\structA \mid x\in \complexB)$ with the same signature as $\structB$ consider first if the sets $\Relationsymbol_x^\structA$ are disjoint. Otherwise, there is no homomorphism $\structA\to \structB$.
	If they are disjoint, define a complex $\complexA$ with a partially defined map $\partialmap\colon \subsetAa \to \complexB$ as follows: The set of vertices of $\complexA$ is the base set of $\structA$. The set $\subsetAa$ is the union $\bigcup_{x\in \complexB} \Relationsymbol_x^\structA$ and the map $\partialmap$ is defined as $\partialmap(y)=x \iff y\in \Relationsymbol_x^\structA$. Moreover, define a subset of $\complexA$ to be a face, if there is a face relation containing a tuple which contains all elements of this set. This way, every relation induces at most $2^{\dim(\complexB)+2}$ faces, which is a constant factor.
	
	Given $(\complexA,\subsetAa,\partialmap)$, define the relational structure $\structA=(\setstructA; \facerel_k^\structA \mid k\le \dim(\complexB)+1, \Relationsymbol_x^\structA \mid x\in \complexB)$ as follows: The base set of $\structA$ is the set of vertices of $\complexA$. The relation $\facerel_k$ contains all ordered $(k+1)$-tuples that are contained in a face of $\complexA$. So every face induced at most $(\dim(\complexB)+2)!$ many tuples, which is a constant factor. Furthermore, define $\Relationsymbol_x^\structA$ as $\{y\in \subsetAa\mid \partialmap(y)=x\}$ for all $x\in \complexB$.
	
	The decision problems are equal because a map $\complexA\to \complexB$ is a homomorphism if it preserves all faces of dimension at most $\dim(\complexB)+1$ by Lemma~\ref{LemmaDimensionForHomo}. 
	As the translation of the inputs can be made in polynomial time in both directions, we get the equivalence in complexity.
\end{proof}
\begin{cor} \label{CorollaryKomplexitaet}
	Let $\complexB$ be a finite simplicial complex.
	If it has an idempotent simplicial Taylor polymorphism, then its precolored CSP is in P. Otherwise, it is NP-complete.
\end{cor}
\begin{proof}
	This is a combination of Theorem~\ref{TheoremZhuk}, Lemma~\ref{LemmaPReduktion} and Lemma~\ref{LemmaPolErhalten}.
\end{proof}

\section{Hom-Complex for Structures}
The notion of the simplicial homomorphism complex generalizes for structures. We give the definition and the connection to pp-interpretations.
\subsection{Defining \universal{} structures}

Note that the homomorphisms between two structures again define an abstract simplicial complex:
\begin{defi}
	For two $\signature$-structures $\structA$ and $\structB$ and a subset $\subsetA\subset \setstructA$, define
	the abstract simplicial complex $\Homsc_\subsetA(\structA,\structB)$ as follows: Its vertex set is the set of all maps 
	$\subsetA\to \setstructB$ which extend to a morphism $\structA\to \structB$. A subset $\face$ of this set is a face, if for all maps $i\colon \subsetA \to \face, a\mapsto i_a$, the map $\subsetA \to \structB, a \mapsto i_a(a)$ is also in the vertex set. 
	%
\end{defi}
Note that this definition is designed such that $\structB$ primitively positivly interprets the finite idempotent realization of $\Homsc_\subsetA(\structA,\structB)$.

With this definitions, there is again a notion of restrictions of homomorphisms and of \universal{} and \contractible{}:
\begin{defi} \label{DefinitionUniversalForStructure}
	%
	A $\signature$-structure $\structB$ is \de{\contractible{}}, if for every finite $\signature$-structure $\structA$ and
	every subsets $\subsetA\subset \setstructA$, 
	the complex $\Homsc_{\subsetA}(\structA,\structB)$ has only contractible components.
	
	We call a $\signature$-structure $\structB$ \de{\universal}, if for every finite abstract simplicial complex $\complexC$, there is a finite $\signature$-structure $\structA$ and a subset $\subsetA
	\subset \setstructA$, 
	such that the complex $\Homsc_{\subsetA}(\structA,\structB)$ is homotopy equivalent 
	to $\complexC$.
\end{defi}

This definitions are well-behaved with the similar definitions for simplicial complexes:
\begin{lemma} \label{LemmaUniversalErhalten}
	Let $(\complexB,\facesB)$ be a finite simplicial complex and let $\structB$ be its finite idempotent realization. Then, 
	\begin{itemize}
		\item $\complexB$ is \contractible{} if and only if $\structB$ is \contractible{} and
		\item $\complexB$ is \universal{} if and only if $\structB$ is \universal{}.
	\end{itemize}
\end{lemma}
\begin{proof}
	The lemma follows directly from the definitions.
\end{proof}

Now, we can add a third condition to Theorem~\ref{TheoreminvPolA}:
\begin{thm} \label{TheoreminvPolB}
	Let $\structB$ be a $\signature$-structure and $\subsetB \subset \setstructB^d$ be a subset of a power. Then, the following are equivalent:
	\begin{enumerate}
		\item \label{NrTheoreminvPolB1}
		There is a primitive positive $\signature$-formula $\phi$ with $d$ many free variables such that for all $x_1,\dots,x_d$ in $\setstructB$, we get $\phi(x_1,\dots,x_d) \iff (x_1,\dots,x_d)\in \subsetB$.
		\item The set $\subsetB$ is invariant under polymorphisms.
		\item \label{NrTheoreminvPolB3}
		There is a second $\signature$-structure $\structA$ and a subset $\subsetA\subset \setstructA$ and a bijection $b\colon \{1,2,\dots,d\} \to \subsetA$ such that the map 
		$$\Homsc_\subsetA(\structA,\structB) \to \setstructB^d, f\mapsto(f(b(1)),f(b(2)),\dots,f(b(d)))$$
		is a bijection from the vertices onto $\subsetB$.
	\end{enumerate}
\end{thm}
\begin{proof}
	It is left to show the equivalence of \ref{NrTheoreminvPolB1} and \ref{NrTheoreminvPolB3}. This translation between homomorphisms and primitive positive formulas is explained in \cite[Section 2.4]{PolyUse}.
\end{proof}
\begin{lemma}\label{LemmaPPIntHomotopie}
	Let $\structA$ and $\structB$ be $\signature$-structures. Let $\subsetA \subset \setstructA$ be a subset. Let $\structD$ be a $\tau$ structure such that $\structD$ pp-interprets $\structB$ with a 1-dimensional interpretation. Then, there exist a $\tau$-structure $\structC$ and a subset $\subsetC \subset \setstructC$ such that $\Homsc_\subsetA(\structA,\structB)$ is homotopy equivalent to $\Homsc_\subsetC(\structC,\structD)$.
\end{lemma}
\begin{proof}
	We fix any 1-dimensional pp-interpretation $\{\psi, \phi_=,\phi_R \mid R\in \signature\}$. Denote with $D'$ the pp-defined subset $\{e\mid \psi(e)\}$ of $\structD$. Recall that there is a surjective map $h\colon D'\to \structB$.
	
	By Theorem~\ref{TheoreminvPolB}, there is a pp-$\tau$-formula $\theta$ with $|\subsetC|$ many free variables representing $\Homsc_\subsetA(\structA,\structB)$. Let $\phi_\theta$ be the $\sigma$-formula, which we obtain from $\theta$ by replacing each relation $R$ with the formula $\phi_R$. Choose $\structC$ and $\subsetC$ such that $\Homsc_\subsetC(\structC, \structD)$ corresponds to $\phi_\theta$ as in Theorem~\ref{TheoreminvPolB}.
	
	Since every free variable of $\phi_\theta$ represents a free variable of $\theta$, the surjection $h\colon D'\to \structB$ induces a surjection $f_1\colon \{\yvect\in (D')^{|\subsetC|}\mid \phi_\theta(\yvect)\} \to \{b\in \structB^{|\subsetC|}\mid \theta(b)\}$ from the solutions of $\phi_\theta$ to the solutions of $\theta$. 
	Choose for every element of $\structB$ a preimage in $D'\subset \structD$. This defines an injection $f_2$ of the solutions of $\theta$ into the solutions of $\phi_\theta$. 
	Both maps induce simplicial maps between the simplicial complexes $\Homsc_\subsetC(\structC,\structD)$ and $\Homsc_\subsetA(\structA,\structB)$ by precomposition. 
	Moreover, $f_1f_2$ is the identity on $\Homsc_\subsetA(\structA,\structB)$ and $f_2f_1$ is homotopic to the identity on $\Homsc_\subsetA(\structA,\structB)$ where the homotopy is given by
	\begin{align*}
		\georeal{\Homsc_\subsetC(\structC,\structD)} \times [0,1] &\to \georeal{\Homsc_\subsetC(\structC,\structD)}\\
		(x,\lambda) &\mapsto \lambda x + (1-\lambda)\georeal{f_1\circ f_2}(x)
	\end{align*}
	Therefore, we get the homotopy equivalence. Note that the line from a point $x$ in the simplicial complex to its image $f_1(f_2(x))$ is always in the complex, because all vectors $\yvect$ in $D'$ we have to analyze this way have the same image in $\structB$ and thus are all in the solution space.
\end{proof}

\subsection{Determine \universal{} Structures}
We can now determine the \universal{} structures.
\begin{lemma}[{\cite[Theorem 2.5]{BooleanCase}, \cite[Lemma 16]{BooleanCaseLIPIcs}}]
    \label{Lemma3SATUniversell}
    $\dSAT$ is \universal{}.
\end{lemma}
\begin{rem*}
    The result in \cite[Theorem 2.5]{BooleanCase} is actually slightly stronger: It is shown that in $\Homsc_{\subsetA}(\structA,\dSAT)$, one can always choose $\subsetA=\setstructA$.
    Note also that the definition of the topological solution space in \cite{BooleanCase} is different, but homotopy equivalent, compared to the definition in this text. This is because they work with the product of topological spaces rather than the product of simplicial complexes. That way, it is possible in \cite{BooleanCase} to get universality up to homeomorphism rather than homotopy for some structures on a binary domain. 
    However, this is not possible in general, with either definition, as we show in Example~\ref{ExampleHomotopyNeeded}, so we can stay with our definition which is a pp-interpretation.
\end{rem*}

\begin{lemma} \label{LemmaConstuctToUniversal}
    Let $\structB$ be a finite $\signature$-structure. Assume that $\structB$ pp-interprets $\dSAT$ in one dimension. Then, $\structB$ is \universal{}.
\end{lemma}
\begin{proof}
    By Lemma~\ref{Lemma3SATUniversell}, $\dSAT$ is \universal{} and by Lemma~\ref{LemmaPPIntHomotopie}, also $\structB$ is.
    %
    %
\end{proof}

\begin{cor} \label{CorollaryNoTaylorUniversal}
    Let $\structA$ be a finite idempotent $\signature$-structure that has no Taylor polymorphism.
    Then, $\structA$ is \universal{}.
\end{cor}
\begin{proof}
    If $\structA$ has no Taylor polymorphism, then there is a clone homomorphism $\Pol(\structA) \to \Proj$ by Theorem~\ref{TheoremCyclicTerms}. By Lemma~\ref{LemmaConsturctInterpret}, this implies that $\structA$ pp-interprets $\dSAT$ in one dimension. Finally, $\structA$ is \universal{} by Lemma~\ref{LemmaConstuctToUniversal}.
\end{proof}
\section{Conclusion}
We have now all results together to show the Theorems from Section~\ref{SectionMainTheorem}.
We start with Theorem~\ref{TheoremExtra2}:
\begin{thm}[Theorem~\ref{TheoremExtra2}]
	The following are equivalent for a simplicial complex $\complexB$.
	\begin{enumerate}
		\item \label{NrTheoremExtra2-A}
		The complex $\complexB$ is \universal{}.
		\item \label{NrTheoremExtra2-B}
		The complex $\complexB$ is not \contractible{}.
		\item \label{NrTheoremExtra2-C}
		The complex $\complexB$ has an idempotent Taylor polymorphism.
		\item \label{NrTheoremExtra2-D}
		There is an idempotent Siggers polymorphism \cite{Siggers}, that is a simplicial homomorphism $\siggers\colon \complexB^6\to \complexB$ such that 
		\begin{align*}
			\forall x \in \complexB: &&\siggers(x,x,x,x,x,x)&=x \text{ and}\\
			\forall x,y,z \in \complexB: && \siggers(x,x,y,y,z,z)&=\siggers(z,y,x,z,y,x)
		\end{align*}
		hold.
		\item \label{NrTheoremExtra2-E}
		There are idempotent cyclic Taylor polymorphism, that are simplicial homomorphisms $\cyclic\colon \complexB^n\to \complexB$  for all prime numbers $n>|\complexB|$ such that 
		\begin{align*}
			\forall x \in \complexB: &&\cyclic(x,\dots ,x)&=x \text{ and}\\
			\forall x_0,x_1,...,x_{n-1} \in \complexB: && \cyclic(x_0,x_1,\dots,x_{n-1})&=\cyclic(x_1,\dots,x_{n-1},x_0)
		\end{align*}
		hold.
	\end{enumerate}
\end{thm}
\begin{proof} We prove the following implications:
	\begin{description}
		\item[\ref{NrTheoremExtra2-C} $\iff$ \ref{NrTheoremExtra2-D} $\iff$ \ref{NrTheoremExtra2-E}] 
		This was shown in Theorem~\ref{TheoremCyclicTerms}.
		\item [\ref{NrTheoremExtra2-E} $\implies$ \ref{NrTheoremExtra2-A}] 
		Choose any $(\complexA, \subsetA,\subsetAa,\partialmap)$. If $\complexB$ has a cyclic Taylor polymorphism then $\Homsc_{\subsetA, \partialmap}(\complexA,\complexB)$ has one as well by Corollary~\ref{KorollarClonhomozuHom}. Thus, it is homotopy equivalent to a discrete set by Theorem~\ref{TheoremSimplicialCyclic}. As $(\complexA, \subsetA,\subsetAa,\partialmap)$ were arbitrary, $\complexB$ is \contractible{}.
		\item [\ref{NrTheoremExtra2-A} $\implies$ \ref{NrTheoremExtra2-B}] 
		Since $\complexB$ is universal, there are $(\complexA, \subsetA,\subsetAa,\partialmap)$ such that $\Homsc_{\subsetA, \partialmap}(\complexA,\complexB)$ is homotopy equivalent to the sphere $\IS^1$. As this is not contractible, $\complexB$ is not \contractible{}.
		\item [\ref{NrTheoremExtra2-B} $\implies$ \ref{NrTheoremExtra2-C}] 
		We prove this by contraposition. Let $\structB$ be the finite idempotent realization of $\complexB$. Then, $\structB$ has no Taylor polymorphism by Lemma~\ref{LemmaPolErhalten}. Thus, $\structB$ is \universal{} by Corollary~\ref{CorollaryNoTaylorUniversal}. Finally, $\complexB$ is \universal{} by Lemma~\ref{LemmaUniversalErhalten}.
		\qedhere
	\end{description}
\end{proof}
The other Theorems are now direct consequences of the equivalence:
\begin{proofof}{Theorem~\ref{TheoremMain}}
	This is the equivalence \ref{NrTheoremExtra2-A} $\iff$ \ref{NrTheoremExtra2-B}.
\end{proofof}

\begin{proofof}{Theorem~\ref{TheoremExtra}}
	The criteria in Theorem~\ref{TheoremExtra} follow similar to the implication \ref{NrTheoremExtra2-E} $\implies$ \ref{NrTheoremExtra2-A}. We do a contraposition: Assume that $\complexB$ is not \universal{}. Then, it is \contractible{} and has cyclic Taylor polymorphisms. So, for every $(\complexA,\subsetA,\partialmap)$, also $\Hom_{\subsetA, \partialmap}(\complexA,\complexB)$ has cyclic Taylor polymorphisms by Corollary~\ref{KorollarClonhomozuHom}. Thus, $\complexB$ and $\Hom_{\subsetA, \partialmap}(\complexA,\complexB)$ are homotopy equivalent to a discrete set by Theorem~\ref{TheoremSimplicialCyclic}.
\end{proofof}

\begin{proofof}{Theorem~\ref{TheoremExtra3}}
	This is a consequence of \ref{NrTheoremExtra2-A} $\iff$ \ref{NrTheoremExtra2-C} and Corollary \ref{CorollaryKomplexitaet}.
\end{proofof}

\section{Generalizations and Applications}

We give a short outlook on possible generalizations of Theorem~\ref{TheoremMain}. While we often get only partial results, this also includes a dichotomy for idempotent relational structures, Theorem~\ref{TheoremDualityForStructures}, and counterexamples to some tempting but wrong generalizations.


\subsection{Omitting the Assumption of Idempotence}

Whenever we looked on polymorphisms of simplicial complexes, we always looked on the idempotent polymorphisms instead of the set of all polymorphisms.

The reason behind restricting the set of all polymorphisms of a simplicial complex is, that every simplicial complex has by definition constant polymorphisms on every vertex, so the full polymorphism clone has no information. In the sections where we looked at polymorphisms of structures, this is not needed and thus we look at all polymorphisms.

There are other possibilities to restrict the polymorphisms of simplicial complexes: In the paper \cite{TopologyAndAdjunction}, a group action on the complex was introduced and the view restricted to polymorphisms which preserve this action. They used that their structures were undirected graphs and flipping the two points of the graph consisting of a single edge is a nontrivial automorphism.


We are able to use our result to get a fixed-point Theorem in this setting which we will use in Section~\ref{SectionHellNesetril} to reprove the Hell-Nešetřil-Theorem. 
To prove it, we need a generalization of Brouwer's fixed-point theorem, which we give with a proof by Borsuk \cite{Borsuk}.
\begin{thm}[A generalization of Brouwer's fixed-point theorem] \label{TheoremBrouwer}
	Let $\complexA$ be a contractible simplicial complex and $\auto\colon \georeal{\complexA} \to \georeal{\complexA}$ a continuous automorphism. Then, $\auto$ has a fixed-point.
\end{thm}
\begin{proof}
	Every contractible simplicial complex is the retract of a simplex. By Brouwer's fixed-point theorem for convex sets \cite{Brouwer}, every automorphism on a simplex has a fixed point. This property is stable under taking a retract \cite[Theorem 7]{Borsuk}.
\end{proof}

This now allows us to conclude a statement about abstract simplicial complexes and polymorphisms.
\begin{thm} \label{TheoremWithAutomorphism}
	Let $(\complexA,\facesA)$ be a finite abstract simplicial complex and $\auto$ a simplicial automorphism of it. Assume that both
	\begin{itemize}
		\item $\Pol(\complexA,\facesA, \auto)$, the set of all simplicial polymorphisms of $\complexA$ that commute with $\auto$ (but do not need to be idempotent), contains no Taylor polymorphism and
		\item there is $x\in \complexA$ such that $x$ and $\auto(x)$ are in the same connected component.
	\end{itemize}
	Then, there is a nonempty face $\face\in \facesA$ that is invariant under $\auto$.
\end{thm}
\begin{proof}
	With endomorphisms and polymorphisms of $(\complexA,\facesA, \auto)$ we denote simplicial maps that commute with $\auto$.
	Let $\complexB\subset \complexA$ be a core of $(\complexA,\facesA, \auto)$.
	That is by definition a subset of $\complexA$ such that there is a surjective endomorphism $f$ from $(\complexA,\facesA, \auto)$ to $(\complexB,\facesA|_{\complexB}, \auto|_{\complexB})$ and every endomorphism of $(\complexB,\facesA|_{\complexB}, \auto|_{\complexB})$ is an automorphism. Such a core exists as $\complexA$ is finite. 
	
	It is well known that the clones $\Pol(\complexA,\facesA, \auto)$ and $\Polidem(\complexB,\facesA|_{\complexB}, \auto|_{\complexB})$, the idempotent part of the core, are homomorphically equivalent (see for example \cite[Theorem 16-17]{PolyUse}) and thus we get in this case that $(\complexB,\facesA|_{\complexB})$ has an idempotent Taylor polymorphism.
	
	Consider the connected component of $f(x)$ in $\complexB$. Since $\auto(x)$ and $x$ are in the same connected component, also $f(\auto(x))=\auto(f(x))$ is in this component. Thus, it is invariant under $\auto$. It is moreover contractible by Theorem~\ref{TheoremTaylorContractible}. Recall that we can extend every simplicial map to the geometric realization, Definition~\ref{DefinitionGeometricRealization}. By extending $\auto$, we get an automorphism of the geometric realization of $\complexB$, which respects the contractible, compact connected component of $f(x)$. By Theorem~\ref{TheoremBrouwer}, it has a fixed point $y$.
	
	Let $\face$ be the smallest face containing the fixed point. Then, $y$ is given by a convex combination of all vertices in $\face$ with nonzero coefficients. As this description of $y$ is unique, $\auto$ has to permute the summands and thus the vertices of $\face$. In conclusion, $\face$ is invariant of $\auto$.
\end{proof}
\subsection{Dichotomy for Relational Structures}  
\label{SectionRelationalStructuresDichotomy}

It is widely studied and there are many known equivalent conditions on when a structure has a Taylor polymorphism. However, it is still an open research problem to find more criteria as this dividing line is highly relevant in complexity theory.
Let us briefly recall some equivalent statements. None of the next Theorem is new and there are more equivalent properties, see e.g. \cite{PolyUse}.


\begin{thm}
For a finite structure $\structA$ the following is equivalent:
\begin{enumerate}
    \item \label{TheoremEquivalentTaylorNrTaylor}
    $\structA$ has a Taylor polymorphism.
    \item \label{TheoremEquivalentTaylorNrInterpretK3}
    $\structA$ cannot pp-interpret $K_3$, the complete graph on three vertices without loops.
    \item \label{TheoremEquivalentTaylorNrInterpret3Sat}
    $\structA$ cannot pp-interpret 3-SAT.
    \item \label{TheoremEquivalentTaylorNrInterpretAll}
    There is a finite structure, which $\structA$ cannot pp-interpret.
    \item \label{TheoremEquivalentTaylorNrClonhomo}
    $\Pol(\structA)$ has no a clone homomorphism to $\proj$.
    \item \label{TheoremEquivalentTaylorNrSiggers6}
    $\structA$ has a 6-ary polymorphism $S\colon \structA^6\to \structA$ such that $\forall x,y,z\in \structA: S(x,x,y,y,z,z)=S(z,y,x,z,y,x)$ (6-ary Siggers Polymorphism).
    \item \label{TheoremEquivalentTaylorNrSiggers4}
    $\structA$ has a 4-ary polymorphism $s\colon \structA^4\to \structA$ such that $\forall x,y,z\in \structA: s(x,y,z,z)=s(z,x,x,y)$  (4-ary Siggers Polymorphism).
    \item \label{TheoremEquivalentTaylorNrCyclicP}
    There is a number $n_0$ and for each prime $p$ larger than $n_0$, there is a cyclic polymorphism $P$ of arity $p$, that is a polymorphism such that 
    $$\forall x_0,...,x_{p-1}\in \structA:P(x_0,x_1...,x_{p-1})=P(x_1,...,x_{p-1},x_0)$$
    holds.
    \item \label{TheoremEquivalentTaylorNrCyclic}
    There is a number $n$ and a cyclic polymorphism $P$ of arity $n$, that is a polymorphism such that 
    $$\forall x_0,...,x_{n-1}\in A:P(x_0,x_1...,x_{n-1})=P(x_1,...,x_{n-1},x_0)$$
    holds.
\end{enumerate}

If $\structA$ satisfies any of the equivalent conditions, $\CSP(\structA)$ is in P. Otherwise, it is NP-complete.
\end{thm}
\begin{proof}[References for Proofs]
    The Implications $\ref{TheoremEquivalentTaylorNrInterpret3Sat}, \ref{TheoremEquivalentTaylorNrInterpretK3}\implies\ref{TheoremEquivalentTaylorNrInterpretAll}$ are obvious.

    The equivalence of \ref{TheoremEquivalentTaylorNrClonhomo} and \ref{TheoremEquivalentTaylorNrInterpret3Sat} is a part of Lemma~\ref{LemmaConsturctInterpret}. The equivalence to \ref{TheoremEquivalentTaylorNrInterpretAll} and \ref{TheoremEquivalentTaylorNrInterpretK3} follows similar using that there is a clone homomorphism $\Pol(K_3)\to \proj$ \cite[Example 34]{PolyUse} and a clone homomorphism from $\proj$ into every clone.

    The equivalence of  \ref{TheoremEquivalentTaylorNrTaylor}, \ref{TheoremEquivalentTaylorNrClonhomo}, \ref{TheoremEquivalentTaylorNrSiggers6}, \ref{TheoremEquivalentTaylorNrSiggers4},  \ref{TheoremEquivalentTaylorNrCyclicP} and \ref{TheoremEquivalentTaylorNrCyclic} follows from Theorem~\ref{TheoremCyclicTerms}.
    
    The complexity theoretic characterization is Theorem~\ref{TheoremZhuk}.
\end{proof}

Recall that we also defined \universal{} and \contractible{} for structures in Definition~\ref{DefinitionUniversalForStructure} and they coincide with the definition for simplicial complexes as described in Lemma~\ref{LemmaUniversalErhalten}.


Before giving the topological dichotomy, consider these two implications:
\begin{thm} \label{TheoremStructuresClassificationImplications}
    Let $\structA$ be a finite idempotent structure.
    \begin{enumerate}
        \item If $\structB$ is a core (i.e. every endomorphism is an automorphism) and has a Taylor polymorphism, then it is \contractible{}.
        \item If $\structB$ is idempotent and has no Taylor polymorphism, then it is \universal{}.
    \end{enumerate}
\end{thm}
\begin{proof}
	For the first direction, note that every core with a Taylor polymorphism also has an idempotent Taylor polymorphism. Consider any complex of the form $\Homsc_\subsetA(\structA,\structB)$. 
	By the clone homomorphism, $\Homsc_\subsetA(\structA,\structB)$ has an idempotent Taylor polymorphism. By Theorem~\ref{TheoremTaylorContractible}, it is contractible.

    The other direction is Corollary~\ref{CorollaryNoTaylorUniversal}.
\end{proof}

This gives a topological duality theorem on finite idempotent structures:
\begin{thm} \label{TheoremDualityForStructures}
	Let $\structA$ be a finite idempotent structure. Then either
	\begin{enumerate}
		\item $\structB$ is \contractible{}. In this case, it has a Taylor polymorphism and its CSP is in P. Or
		\item $\structB$ is \universal{}. In this case, it has no Taylor polymorphism and its CSP is NP-complete.
	\end{enumerate}
\end{thm}
\begin{proof}
	This is a combination of Theorem~\ref{TheoremStructuresClassificationImplications} and~\ref{TheoremZhuk}.
\end{proof}

\begin{rem}
	The assumptions core respectively idempotent in Theorem~\ref{TheoremStructuresClassificationImplications} are necessary, see Examples~\ref{ExampleIdempotentNeeded} and~\ref{ExampleCoreNeeded} below. However, that is no big constraint as every CSP is polynomial time equivalent to the CSP of an idempotent structure \cite[Theorem 16-17]{PolyUse}.
\end{rem}

Every idempotent structure without a Taylor polymorphism is \universal{}. The assumption idempotent is necessary:
\begin{example}
	\label{ExampleIdempotentNeeded}
	The structure on the set $\{0,1\}$ with the 3-ary relation not all equal respectively $\{0,1\}^3\setminus \{(0,0,0),(1,1,1)\}$ is NP-complete. However it is not \universal{} as it cannot represent three contractible components.
\end{example}
\begin{proof}
	The idea is that every complex represented by the structure has an automorphism by permuting $\{0,1\}$. A complex with an odd number of connected components thus has a fixed component and by Theorem~\ref{TheoremBrouwer} a fixed face. The only faces that are invariant are faces that contain two opposite points of $\{0,1\}^n$. But these faces exist only if all points in $\{0,1\}^n$ exist and thus are connected.
\end{proof}

Every core structure with a Taylor polymorphism is \contractible{}. The assumption core is necessary:
\begin{example}
	\label{ExampleCoreNeeded}
	Take on the set $\{0,1\}$ the controlled 3-SAT structure, that are the relations 
	\begin{align*}
		&\emptyset \text{ (the empty relation)}, \\
		N\coloneqq {}&\{(0,1,0),(0,0,1),(1,1,1)\} \text{ and}\\
		&\{0,1\}^3 \setminus \{(0,0,0)\}
	\end{align*}
	This structure is \universal{} (and thus not \contractible{}), but has a constant polymorphism by mapping everything to 1.
	
	Moreover, the same structure without the empty relation is neither \universal{} nor \contractible{}. 
\end{example}
\begin{proof}
	Every nonempty simplicial complex can be obtained from this structure without the empty relation up to homotopy by representing it first by $\dSAT$ in a way such that the constant 1-map is in the complex and every variable is contained in a $\NOT$-relation. Now, we add one point $c$ to the instance of $\dSAT$ (but not to the subset $\subsetA$) and replace every occurrence of $\NOT(x,y)$ by $N(c,x,y)$. The empty complex can be obtained by the empty relation. Without the empty relation, the empty complex cannot be defined.
\end{proof}

\subsection{Omitting the Homotopy}
In the definition of \universal{}, we only show that we can obtain every complex up to homotopy and not up to a homeomorphism. However, there are \universal{} structures, for which we cannot get the circle $\IS^1$ up to homeomorphism.
\begin{example}[The need of Homotopies] \label{ExampleHomotopyNeeded}
    The structure $\structB$ on $\{1,2,3,4\}$ with the unary relations $\{1\}, \{2\}, \{3\}, \{4\}$ and the higher relations 
    \begin{align*}
        \text{wide NOT} &\coloneqq \{(x,y)\mid x\in \{1,2\}\iff y\in \{3,4\}\} \text{ and} \\
        \text{wide OR} &\coloneqq \{(x,y,z)\mid (x,y,z)\notin \{1,2\}^3\}\}
    \end{align*}
    is idempotent and \universal{} but $\Homsc_{\subsetA}(\structA,\structB)$ is never homeomorphic to the sphere.
\end{example}
\begin{proof}
    It is \universal{} as it pp-interprets $\dSAT$. Take any $(\structA,\subsetA)$, such that $\Homsc_{\subsetA}(\structA,\structB)$ is homotopy equivalent to the sphere. Every point in $\subsetA$ is either fixed, in which case we can remove it from $\subsetA$ and get a homeomorphic complex. Or the point only occurs in the wide OR and wide NOT relation. This happens for at least one point. Then, it has always at two possible values, namely $1$ and $2$ or $3$ and $4$ (or all four). But this implies that $\Homsc_{\subsetA}(\structA,\structB)$ has $[0,1]$ as a direct factor and thus is not the sphere.
\end{proof}

We can get the same problem with a simplicial complex:
\begin{example}
	Consider the simplicial complex with vertex set $\{a_1,a_2,b_1,b_2,c_1,c_2\}$ and as faces the sets $\{\{a_1,a_2,b_1,b_2\},\{b_1,b_2,c_1,c_2\},\{a_1,a_2,c_1,c_2\}\}$ and all of its subsets.
	$$
	
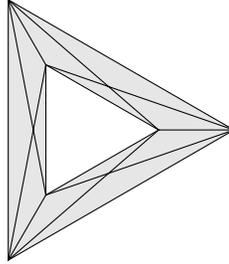
\begin{tikzpicture}[even odd rule]
		\filldraw[fill=gray!20!white] 
		(0:2)--(120:2)--(240:2)--cycle
		(0:1)--(120:1)--(240:1)--cycle
		;
		\draw (0:2)--(120:1)--(240:2)--(0:1)--(120:2)--(240:1)--cycle
		(0:1)--(0:2)
		(120:1)--(120:2)
		(240:1)--(240:2)
		;
	\end{tikzpicture}
	$$
	This complex is \universal{} as it is not contractible, but $\Homsc_{\subsetA}(\structA,\structB)$ is never homeomorphic to the sphere.
\end{example}
The proof is similar to Example~\ref{ExampleHomotopyNeeded}.

\subsection{Omitting Projections}
In our notion of \universal{}, we include the possibility that $\subsetA \subsetneq \structA$ is a proper subset. This is needed to show that every structure without Taylor polymorphism is \universal{}.
\begin{example}[The need of Projections, by {\cite[Lemma 31]{BooleanCaseLIPIcs}}]
    Take the structure 1-in-3-sat on the set $\{0,1\}$ consisting of the relation
    $$
        \{(0,0,1), (0,1,0),(1,0,0)\}.
    $$
    This is \universal{}, but for $\subsetA=\structA$, we get that every connected component of $\Homsc_{\subsetA}(\structA, \text{1-in-3-sat})$ is contractible.
\end{example}
For simplicial complexes, it might be possible to omit this assumption.

\subsection{Converse to Theorem~\ref{TheoremTaylorContractible}}
By Theorem~\ref{TheoremTaylorContractible}, a simplicial complex with a Taylor polymorphism is homotopy equivalent to a discrete set. The converse is false:
\begin{example}
    The 2-dimensional simplicial complex
    $$
        
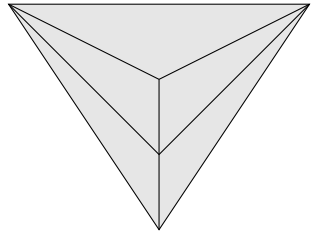
\begin{tikzpicture}
            \draw[fill=gray!20!white] (0,0) -- (2,3) -- (-2,3) -- cycle;
            \draw (0,0) -- (0,1) -- (0,2);
            \draw (2,3) -- (0,1) -- (-2,3);
            \draw (2,3) -- (0,2) -- (-2,3);
        \end{tikzpicture}
    $$
    is contractible but it is not Taylor nor \contractible{}.
\end{example}
\begin{proof}
    Consider the full subsimplex of all vertices that are directly connected to the lowest vertex. This subsimplex is not contractible and has by Theorem~\ref{TheoremTaylorContractible} no Taylor polymorphism. As the above Simplex pp-interprets the subsimplex in one dimension, it also has no Taylor polymorphism and is therefore \universal{} and not \contractible{}.
\end{proof}
Note that we have a converse of Theorem~\ref{TheoremTaylorContractible} in the sense that if $\complexB$ has no Taylor polymorphism then there is a non-contractible complex of the form $\Homsc_{\subsetA,\partialmap}(\complexA,\complexB)$ by Theorem~\ref{TheoremMain}.

\subsection{Infinite Complexes}

There are two interesting generalizations of finite simplicial complexes: Infinite simplicial complexes and pairs of finite simplicial complexes, similar to the generalizations of finite constraint satisfaction problems to infinite constraint satisfaction problems and promise constraint satisfaction problems (PCSPs).

For the infinite case, the initial question is obvious:
\begin{question}
    Does a version of Theorem~\ref{TheoremMain} holds for infinite structures?
\end{question}
A key point to the proof is Theorem~\ref{TheoremTaylorContractible}. Even tough we have two proofs for it, we use the assumption finite in both proofs. For the proof written in this paper, it is used in the cyclic terms Theorem.
However, there is a notable case in which Theorem~\ref{TheoremTaylorContractible} generalizes:
\begin{defi}
	A simplicial complex is called \de{locally finite} if every vertex is contained in only finitely many faces.
\end{defi}
\begin{thm}
	Let $\complexA$ be a locally finite abstract simplicial complex that has a simplicial idempotent Taylor polymorphism. Then, every connected component of $A$ is contractible.
\end{thm}
\begin{proof}
	Note that for any point $x\in \complexA$ and any $n\in \IN$, the set of all points with distance at most $n$ from $x$ is invariant under polymorphisms. It is moreover finite as $\complexA$ is locally finite and thus by Theorem~\ref{TheoremTaylorContractible} contractible. To show that also every component of infinite size is contractible, it suffices by Whitehead's theorem to show that every map from a sphere is contractible. However, this is true as the image of a sphere can be approximated by a finite connected subcomplex and is thus contained in a finite, contractible set.
\end{proof}

For the other direction, that every complex with a Taylor polymorphism is \contractible{}, there is hope as there is a topological Birkhoff Theorem for infinite structures \cite[Section 9.5]{theBodirsky}. 

For infinite structures, non-idempotent polymorphisms are of a wider interest. The in some sense simplest infinite structures are oligomorphic structures which are never idempotent but should ideally be included in the theorem.

\subsection{Pairs of Complexes}  

For a pair of structures, the situation is again more complicated. At first, we need to consider the equivalent of ``an idempotent structure" for a pair of complexes to get able to formulate a generalization of Theorem~\ref{TheoremMain}.
\begin{defi}
    We call a triple $(\complexB,\complexC,\phi)$ of two simplicial complexes and a simplicial map $f\colon \complexB \to \complexC$ \de{\contractible{}}, if for all tuples $(\complexA,\subsetA,\subsetAa,\partialmap\colon \subsetAa\to \complexB)$ the post composition with $\phi$, considered as map from $\Homsc_{\subsetA,\partialmap}(\complexA,\complexB)$ to $\Homsc_{\subsetA,\phi \circ \partialmap}(\complexA,\complexC)$, is homotopic to a locally constant map.
\end{defi}
Equivalently, for each $C$, the composition with $\phi$ factors through a set that is homotopy equivalent to a discrete set.
Note that this definition generalizes Definition~\ref{DefinitionContractible} as $(\complexB,\complexB,\id_\complexB)$ is \contractible{} if and only if $\complexB$ is, by definition.

The next problem is, that the set of maps $\complexB^n\to \complexC$ which extend $\phi$ are not a clone, but a minion. A minion is a generalization of a clone. For a precise definition of Minion, promise constraint satisfaction problems and bounded essential arity, we refer the reader to \cite{BBKO}.
\begin{question}
    Consider a triple $(\complexB,\complexC,\phi)$. Does the minion
    $$
        \Pol_\phi(\complexB,\complexC)\coloneqq \{f\colon \complexB^n\to \complexC\mid n \in \IN, f(x,x,\dots,x)=\phi(x) \}
    $$
    defines whether the triple $(\complexB,\complexC,\phi)$ is \contractible{}?
\end{question}
The same question can be asked for general structures where \contractible{} can be defined in an analogous way. A more concrete question is
\begin{question}
    Consider a triple $(\structB,\structC,\phi)$ of two structures and a map. Consider the statements
    \begin{enumerate}
        \item It is \contractible{}.
        \item The complexity of the promise constraint satisfaction problem described by the minion
        $$
        \Pol_\phi(\structB,\structC)\coloneqq \{f\colon \structB^n\to \structC\mid n \in \IN, f(x,x,\dots,x)=\phi(x) \}
        $$
        is in P.
    \end{enumerate}
    Does one of the statements implies the other?
\end{question}
None of the two directions is obvious and a ``no" as answer for at least one of the directions would be no surprise.

With some additional assumptions, this setting has been studied in \cite{TopologyAndAdjunction}. The next Theorem is a variation of \cite[Theorem 1.4]{TopologyAndAdjunction} and can be proven analogously:
\begin{thm} \label{TheoremBoundedEssArity}
    Consider a triple $(\structB,\structC,\phi)$ of two structures and a map.
    Assume that there is $(\structA,\subsetA)$ such that the homotopy groups of $\Homsc_\subsetA(\structA,\structB)$ and $\Homsc_\subsetA(\structA,\structC)$ are both $\IZ$ and the map induced by composition with $\phi$ on the fundamental groups is not the (constant) $0$ map. 
    Then, the minion $\Pol_\phi(\structB,\structC)$ has bounded essential arity and thus the PCSP is NP-hard.
\end{thm}

It should also be mentioned that the existence of a Taylor polymorphism in the promise setting is neither enough to show that the complexity class is P nor to obtain that the map is homotopic to a constant map, so the canonical generalization of Theorem~\ref{TheoremTaylorContractible} is false:
\begin{example}
    Consider $(\cyclecomplex_9,\cyclecomplex_3,\phi)$ with the complexes from Example~\ref{ExampleCyclepath} and the map $m\colon \{1,2,\dots,9\}^3\to \{1,2,3\}, (x,y,z)\mapsto \lfloor x-y+z\rfloor\pmod{3}$ where $\phi(x)=m(x,x,x)$. The map $m$ has the Maltsev property $m(x,x,x)=m(x,y,y)=m(y,y,x)=x$ and is a Taylor polymorphism in $\Pol_\phi(\cyclecomplex_9,\cyclecomplex_3)$. However, the map $\phi$ is not homotopic to a constant map and the Minion $\Pol_\phi(\cyclecomplex_9,\cyclecomplex_3)$ has bounded essential arity by Theorem~\ref{TheoremBoundedEssArity}.
\end{example}

The above example is motivated by the fact that the circle $\IS^1$ is a noncontractible topological space with a continuous Taylor polymorphism. We can also construct examples where the homomorphism is not a homotopy equivalence:
\begin{example}
	\label{ExampleRPPlane}
    Consider the simplicial complex $(\complexA,\facesA)$ given by the dodecahedron model of the real projective plane, that is
    \begin{align*}
        \complexA &=\{1,2,3,4,5,6\} \\
        \facesA &= \{\{123\},\{134\},\{145\},\{156\},\{162\},\{235\},\{346\},\{452\},\{563\},\{624\} \\&\phantom{{}=\{} \text{and all subsets}\}.
    \end{align*}
\begin{figure}
	$$
	\begin{tikzpicture}
		\draw node (a) at   (0,0) {1};
		\draw node (b2) at   (0:1) {2};
		\draw node (b3) at  (72:1) {3};
		\draw node (b4) at (144:1) {4};
		\draw node (b5) at (216:1) {5};
		\draw node (b6) at (288:1) {6};
		\draw node (c5) at  (36:1.62) {5};
		\draw node (c6) at (108:1.62) {6};
		\draw node (c2) at (180:1.62) {2};
		\draw node (c3) at (252:1.62) {3};
		\draw node (c4) at (324:1.62) {4};
		\draw (a)--(b2) (a)--(b3) (a)--(b4) (a)--(b5) (a)--(b6);
		\draw (b2)--(b3)--(b4)--(b5)--(b6)--(b2);
		\draw (b2)--(c5)--(b3)--(c6)--(b4)--(c2)--(b5)--(c3)--(b6)--(c4)--(b2);
	\end{tikzpicture}
	$$
	\caption{The real projective plane from Example~\ref{ExampleRPPlane}. Opposite edges need to be glued together.}
\end{figure}
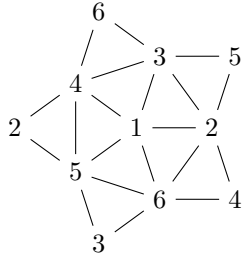
    Let $\complexB=\complexA^3/M$ where $M$ is the equivalence relation identifying $(x,x,y)\sim(x,y,x)\sim (y,x,x)\sim(x,x,x)$ and let $\phi$ be the diagonal map $\phi\colon \complexA\to \complexB$. Then, $\phi$ is trivial in all homotopy groups, but $\phi$ is not contractible as it is nontrivial in $\IF_2$ cohomology \cite{mathOverflow}. Note that $\Pol_\phi(\complexA,\complexB)$ admits a quasi-majority polymorphism $p$ that satisfies $p(x,x,y)=p(x,y,x)=p(y,x,x)=p(x,x,x)$ by mapping $(x,y,z)$ to the equivalence class of $(x,y,z)$.
\end{example}
\section{Applications}

We have already seen some applications in the previous section, including a topological dichotomy for relational structures in Section~\ref{SectionRelationalStructuresDichotomy}. In this section, we want to present two more applications: An alternative proof of the Hell-Nešetřil-Theorem and the connection to social choice functions.

\subsection{Alternative Proof of the Hell-Nešetřil-Theorem}
\label{SectionHellNesetril}
The topological dichotomy gives a possibility to reprove the Hell-Nešetřil-Theorem from graph theory. Here, we use polymorphisms of a simplicial complex, that preserve an automorphism and might be non-idempotent and Theorem~\ref{TheoremWithAutomorphism}. 
A self-contained version of this proof will be published in \cite{OprsalMeyer}.

\begin{thm}[{Hell-Nešetřil-Theorem, \cite[Theorem 1]{Hell-Nesetril}}]
    Let $H$ be an undirected graph without loops.
    If $H$ is bipartite then the $H$-coloring problem is in P.
    If $H$ is not bipartite then the $H$-coloring problem is NP-complete.
\end{thm}
\begin{proof}
    If $H$ is bipartite, then a graph $G$ is $H$-colorable if and only if it is two colorable. This problem is in P.

    If $H$ is non-bipartite, assume that $H$ has a Taylor polymorphism.
    Consider the complex $(\complexA,\facesA)\coloneqq \Homsc(E,H)$, where $E$ is the graph which consists of two vertices $\{v,w\}$ that are linked by an edge. Note that every polymorphism of $H$ induces a polymorphism of $\complexA$. Note moreover that $\complexA$ admits an automorphism $\auto$ which maps the vertex $f\colon E \to H$ to $g\circ f$ where $g\colon E \to E, g(v)\coloneqq w$ and $g(w)\coloneqq v$. This automorphism is invariant under the polymorphisms that are induced by $\Pol(H)$ as the polymorphisms of $H$ act on $\Homsc(E,H)$ from the right while $g$ acts from the left.

    Thus, we get a clone homomorphism from $\Pol(H)$ to $\Pol(\complexA, \facesA, \auto)$. Then, $\Pol(\complexA, \facesA, \auto)$ has a Taylor polymorphism. Thus, by Theorem~\ref{TheoremWithAutomorphism}, there is $\face\in \facesA$ that is invariant under $\auto$. Let $x$ be a vertex in $\face$.
    As $\auto\circ\auto(f)=g\circ g\circ f=f$, we get that either $\auto(\auto(x))=x$ and the $\auto$-orbit of $x$ is $\auto$-invariant and defines a subface of $\face$ and thus a face. Thus, there is an $\auto$-invariant face of size one or two.

    In the first case, we get $\auto(x)=x$ and thus $g\circ x=x$ and therefore $x(v)=x(w)$ so $x$ is a loop in $G$. A contradiction.
    In the second case, we get that $\auto(x)\ne x$, so $x$ is no loop, but $\{x, \auto(x)\}$ is a face. In this case, we get by the definition of the face in $\Homsc(E,H)$ that the maps
    \begin{align*}
        f&=(v\mapsto f(v), w\mapsto f(w))\\
        \auto(f) &= (v\mapsto f(w), w\mapsto f(v))\\
        (v\mapsto f(v), w\mapsto (\auto(f))(w))&= (v\mapsto f(v), w\mapsto f(v))\\
        (v\mapsto (\auto(f))(v), w\mapsto f(w))&= (v\mapsto f(w), w\mapsto f(w))
    \end{align*}
    are in $\Hom(E,H)$. Thus the constant map on $f(v)$ is in $\Hom(E,H)$ and $H$ has a loop on that vertex. A contradiction.

    As we get in both cases a contradiction, our assumption was wrong and $H$ non-bipartite has no Taylor polymorphism. But in this case, the $H$-coloring problem is NP-complete.
\end{proof}
\subsection{Social Choice }

A social choice function is essentially the same as a polymorphism. As the topological dichotomy proven in this paper is closely related to polymorphisms, it is also connected to the existence of social choice functions. This gives some immediate results which we present in this section.

\begin{defi}
	For the purpose of this paper, a \de{choice function on an abstract simplicial complex} $\complexA$ in a set $K$ of agents is a map $f$ from the set of maps $K \to \complexA$ to $\complexA$ satisfying 
	\begin{enumerate}
		\item unanimity, i.e. if the map $K\to\complexA$ is constant $x$ then its $f$-image should be $x$ as well and
		\item the homomorphism property, i.e. if $\subsetA$ is a subset of the maps from $K$ to $\complexA$ such that for each agent $i\in K$, the set $\{g(i)\mid g\in \subsetA\}$ is connected by a face, then the image of $\subsetA$ should also be connected by a face.
	\end{enumerate}
	We call the choice function \de{Taylor} (after \cite{Taylor77}), if for each agent $i$, there exists two sets $S_i$ and $S_i'$ such that $i$ is contained in $S_i$, but not in $S_i'$ and for all $x,y\in \complexA$, the identity
	$$
		f\left(j\mapsto \begin{cases}
			x &\text{if } j\in S_i \\
			y &\text{else}
		\end{cases}\right)
		=
		f\left(j\mapsto \begin{cases}
			x &\text{if } j\in S_i' \\
			y &\text{else}
		\end{cases}\right)
	$$
	holds.
\end{defi}
The Taylor notion essentially claims that the sets $S_i$ and $S_i'$ are treated in a similar way. As $i$ is contained in only one of them, it ensures that $f$ is not the evaluation map on one of the agents or, in the language of social sciences, there is no dictatorship.

Note furthermore that the Taylor notion includes the notions of anonymity, majority and near unanimity from Section~\ref{SectionSocialChoiceIntro}. However, the existence of a Taylor choice function implies by no means the existence of an anonymous or majority choice function. But it is highly related to the topological dichotomy:
\begin{thm}
	An abstract simplicial complex admits a Taylor social choice function if and only if it is \contractible{}. In that case, it is homotopy equivalent to a discrete set. 
\end{thm}
\begin{proof}
	An unanimous social choice function in the set $K$ is equivalent to an idempotent polymorphism of arity $|K|$ by replacing maps $K\to \complexA$ with powers $\complexA^{|K|}$. This association preserves the Taylor property as they are defined in a similar way. Thus the existence of a Taylor choice function is equivalent to the existence of an idempotent Taylor polymorphism which is equivalent to being \contractible{} by Theorem~\ref{TheoremExtra2}. This implies the homotopy equivalence by Theorem~\ref{TheoremExtra}.
\end{proof}



\clearpage

\bibliography{CSP-CW-Quellen} 
\bibliographystyle{alpha}

\end{document}